\documentclass[12pt]{article}
\usepackage{t1enc,amssymb,amsbsy,sal_yor_2}

\begin{document}

\author{Paavo Salminen\\{\small Åbo Akademi,}
\\{\small Mathematical Department,}
\\{\small FIN-20500 Åbo, Finland,} \\
{\small email: phsalmin@abo.fi}
\and
Marc Yor\\
{\small Universit\'e Pierre et Marie Curie,}\\
{\small Laboratoire de Probabilit\'es,}\\
{\small 4, Place Jussieu, Case 188}\\
{\small F-75252 Paris Cedex 05, France}\\
}
\vskip5cm


\title{Perpetual integral functionals as hitting and occupation times}

\maketitle

\begin{abstract}
Let $X$ be a linear diffusion and $f$ a non-negative, Borel
measurable function. We are interested in finding conditions
on $X$ and $f$ which imply that 
the perpetual integral functional
$$
I^X_\infty(f):=\int_0^\infty f(X_t)\, dt
$$
is identical in law with the first hitting time of a point for some other
diffusion. This phenomenon may often be explained using random time
change. Because of some potential
applications in mathematical finance, we are considering mainly the case
when $X$ is a Brownian motion with drift $\mu>0,$ denoted
$\{B^{(\mu)}_t:\ t\geq 0\},$ but it is obvious
that the method presented is more general. We also review the
known examples and give new ones. In particular, results concerning one-sided
functionals
$$
\int_0^\infty f(B^{(\mu)}_t)\,{\bf 1}_{\{B^{(\mu)}_t<0\}} dt\quad
{\rm and}\quad
\int_0^\infty f(B^{(\mu)}_t)\,{\bf 1}_{\{B^{(\mu)}_t>0\}} dt
$$
are presented.

This approach
generalizes the proof, based on the random time change techniques, 
of the fact that the Dufresne functional (this
corresponds to $f(x)=\exp(-2x)),$ playing quite an important
r\^ole in the study of geometric Brownian motion, is identical in law with the
first hitting time for a Bessel process. Another functional
arising naturally in this context is 
$$
\int_0^\infty \big(a+\exp(B^{(\mu)}_t)\big)^{-2}\, dt,
$$
which is seen, in the case $\mu=1/2,$ to be identical in
law with the first hitting time for a Brownian motion with drift
$\mu=a/2.$

The paper is concluded by discussing how the
Feynman-Kac formula can be used to find the distribution of a
perpetual integral functional.
\\ \\
{\rm Keywords:} Time change, Lamperti transformation, Bessel
processes, Ray--Knight theorems, Feynman-Kac formula.
\\ \\ 
{\rm AMS Classification:} 60J65, 60J60, 60J70.
\end{abstract}

\section{Introduction and summary of the results}
\label{sec0}
Let $B^{(\mu)}=\{ B^{(\mu)}_t:=B_t+\mu t \,:\, t\geq
0\}$ be a Brownian motion with drift $\mu>0$ and $f$ a
non-negative measurable function. Encouraged by a number of
examples listed below, we wish to gain better understanding when
an integral functional of the type
$$
I_\infty(f):=\int_0^\infty f(B^{(\mu)}_s)\,ds
$$
is identical in law with the first hitting time of a point for some other
diffusion. Clearly, we can pose an analogous question for an
arbitrary diffusion instead of $B^{(\mu)}.$ In fact, the results
in Section \ref{sec1} are fairly easily extended for arbitrary
transient diffusions determined by a stochastic differential
equation. Our interest in the particular case with $B^{(\mu)}$ is
motivated by the numerous studies and results associated to the
functional
\begin{equation}
\label{du} \int_0^\infty  \exp(-2B_s^{(\mu)})\, ds.
\end{equation}
This functional was first considered by Dufresne in
\cite{dufresne90} where it is seen, among other things, how the
functional (\ref{du}) arises as a perpetuity after a limiting
procedure in a discrete model.

We now review some cases of perpetual integral functionals which are
identical in law with the first hitting time. Let
$$
H_a(Z):=\inf\{t:\ Z_t=a\}
$$
denote the first
hitting time of the point $a$ for a diffusion $Z.$

{\bf 1)}  In Yor \cite{yor92b} (see \cite{yor01} for
an English translation) it is shown that for
the Dufresne functional (\ref{du})
we have
\begin{equation}
\label{d-y}
\int_0^\infty  \exp(-2B_s^{(\mu)})\,
ds
\quad{\mathop=^{\rm{(d)}}}\quad
H_0(R^{(\delta)}),
\end{equation}
where $R^{(\delta)}$ is a Bessel process of dimension
$\delta=2(1-\mu)$ started at 1,
and $\displaystyle{{\mathop=^{\rm{(d)}}}}$ reads
"is identical in law with".
Recall also that
\begin{equation}
\label{d-y1}
\int_0^\infty \exp\bigl(-2B^{(\mu)}_s\bigr)\,ds\quad{\mathop=^{\rm{(d)}}}\quad
{1\over{2\gamma_\mu}},
\end{equation}
where $\gamma_\mu$ is a gamma-distributed random variable
with parameter $\mu.$ We refer to Szabados and Sz\'ekely 
\cite{szabadosszekely03} for a discussion of Dufresne's functional 
for random walks.

{\bf 2)} The Ciesielski--Taylor identity:
\begin{equation}
\label{c-t}
\int_0^\infty {\bf 1}_{\{R^{(\delta+2)}_s< 1\}}\, ds
\quad{\mathop=^{\rm{(d)}}}\quad
H_1(R^{(\delta)}),
\end{equation}
where $R^{(\delta)}$ and $R^{(\delta+2)}$ are Bessel processes of
dimension $\delta>0$ and $\delta+2,$ respectively,
started at 0.
For a proof, see Williams \cite{williams79} p. 159 and 211, and Yor
\cite{yor91}, \cite{yor92c} p. 50; in the case $\delta=1$ there is a pathwise
explanation due to D. Williams. We refer also to Getoor and
Sharpe \cite{getoorsharpe79} p. 98, and
to Biane \cite{biane85}
for a generalization to a vast class of pairs of diffusions.

{\bf 3)} The identity due to Biane \cite{biane85} and Imhof
\cite{imhof86}:
\begin{equation}
\label{i}
\int_0^\infty {\bf 1}_{\{B^{(\mu)}_s< 0\}}\, ds
\quad{\mathop=^{\rm{(d)}}}\quad
H_\lambda(B^{(\mu)}),
\end{equation}
where $\lambda$ is a random variable independent of $B^{(\mu)}$
and exponentially distributed
with parameter $2\mu$ and $B^{(\mu)}_0=0.$

{\bf 4)} The identity due to Donati-Martin and Yor
\cite{donatimartinyor97} p. 1044:
\begin{equation}
\label{dm-y}
\int_0^\infty \frac{ds}{\exp(2R^{(3)}_s)-1}
\quad{\mathop=^{\rm{(d)}}}\quad
H_{\pi/2}(R^{(3)}).
\end{equation}
where $R^{(3)}$ is a three-dimensional Bessel process started from
0. See \cite{donatimartinyor97} also for a probabilistic
explanation of (\ref{dm-y}).

In Section \ref{sec1} of this paper we present a general method
based on It\^o's and Tanaka's formulae and random time change
techniques which connects the distribution of a perpetual integral
functional to the first hitting time. This method gives us the identity
(\ref{d-y}) and also the following (which could be called the
reflecting counterpart of (\ref{d-y})):
\begin{equation}
\label{ps2} 
\int_0^{\infty} \exp(-2a\,B^{(\mu)}_s)\,{\bf
1}_{\{B^{(\mu)}_s> 0\}}\, ds \quad{\mathop=^{\rm{(d)}}}\quad
H_{1/a}(R^{(\delta)}),\quad \delta=2\mu/a.
\end{equation}
where the Bessel process $R^{(\delta)}$ is started at 0 and, in
the case $0<\delta<2,$ reflected at 0.
However, the simplest case emerging from our approach leads us
to the identity
\begin{equation}
\label{y}
\int_0^{\infty}  (a+\exp(B^{(1/2)}_s))^{-2}\, ds
\quad{\mathop=^{\rm{(d)}}}\quad
H_{r}(B^{(a/2)}),
\end{equation}
where $a>0,$ $B^{(a/2)}_0=0$ and $r=\frac{1}{a}\log(1+a).$
The reflecting counterpart of (\ref{y}) is
\begin{equation}
\label{yref}
\int_0^{\infty}
\frac{{\bf
1}_{\{B^{(1/2)}_s> 0\}}}
{(a+\exp(B^{(1/2)}_s))^{2}}\, ds
\quad{\mathop=^{\rm{(d)}}}\quad
H_{r}(\tilde B^{(a/2)}),
\end{equation}
where $a>0,$ $\tilde B^{(a/2)}$ is reflecting Brownian motion with drift $a/2$
started from 0 and $r$ is as above. 
In Section 3.3, when analyzing the functional on the left hand side of 
(\ref{y}), we also find the diffusion with the first hitting time
identical in law with the functional 
$$
\int_0^{\infty}  (a+\exp(B^{(\mu)}_s))^{-2}\, ds,
$$
where $\mu>0$ is arbitrary. The Laplace transform of this functional 
can be expressed in terms of Gauss' hypergeometric functions 
(see \cite{borodinsalminen04}).


We have not investigated or constructed a discrete model (as is done
in \cite{dufresne90} for the functional in (\ref{du})) which would
lead to the
functional in (\ref{y}). 
Notice, however, that
$$
\int_0^{\infty} (1+\exp(B^{(\mu)}_s))^{-2}\, ds
=
\int_0^{\infty} \exp(-2\,B^{(\mu)}_s)
(1+\exp(-B^{(\mu)}_s))^{-2}\, ds
$$
and, hence, this functional can be considered as a modification of
Dufresne's functional such that the discounting is
bounded (we call this modified functional a translated Dufresne's
functional).  In fact, using the results in Salminen and Yor
\cite{salminenyor03}, where the integrability properties of perpetual
functionals are discussed, it is seen that while
Dufresne's functional does not have
moments of order $m\geq\mu$ (cf. (\ref{d-y1})),
which is perhaps unrealistic from an economical point of view,
the functional in (\ref{y}) has some
exponential moments - being in this respect more appropriate.

For functionals restricted to the negative half line we cannot in
general have similar descriptions in terms of first hitting times. A
typical example is the identity (\ref{i}) above. However, the
Lamperti transformation allows us to connect exponential
functionals to the occupation times for Bessel processes. In
Section \ref{ct-revis} we show the identity
\begin{equation}
\label{dref2}
\int_0^{\infty}
\exp(-2a\,B^{(\mu)}_s)\,
\,
{\bf 1}_{\{B^{(\mu)}_s< 0\}}\, ds
\quad{\mathop=^{\rm{(d)}}}\quad
\int_0^{\infty}
{\bf 1}_{\{R^{(\delta_2)}_s> 1/a\}}\, ds,
\end{equation}
where $R^{(\delta_2)}$ is a Bessel process
of dimension $\delta_2=2(1-\frac{\mu}{a})$ started at
$1/a.$ Notice that $R^{(\delta_2)}$ hits 0 in finite
time, and, in case $0<\mu<a,$ we take 0 to be a killing
boundary point.
Further, also by the Lamperti time change, we have
\begin{equation}
\label{dref3}
\int_0^{\infty}
\exp(2a\,B^{(\mu)}_s)\,
\,
{\bf 1}_{\{B^{(\mu)}_s< 0\}}\, ds
\quad{\mathop=^{\rm{(d)}}}\quad
\int_0^{\infty}
{\bf 1}_{\{R^{(\delta_3)}_s< 1/a\}}\, ds,
\end{equation}
where $R^{(\delta_3)}$ is a Bessel process
of dimension $\delta_3=2(1+\frac{\mu}{a})$ started at
$1/a.$ This identity was first observed and proved
in Yor \cite{yor93}
(\cite{yor01} p. 133) by different methods.

Finally, we recall the recent works by the
second author, jointly with H. Matsumoto, see \cite{matsumotoyor99a},
\cite{matsumotoyor00},
\cite{matsumotoyor01},
in which the variable
\begin{equation}
\label{t}
\int_0^t \exp\bigl(-2B^{(\mu)}_s\bigr)\,ds, \quad \mu>0,
\end{equation}
plays an essential role in obtaining a variant
of Pitman's theorem. Indeed, it
is proved that the process
$$
Z^{(\mu)}_t:=\exp\bigl(B^{(\mu)}_t\bigr)
\int_0^t \exp\bigl(-2B^{(\mu)}_s\bigr)\,ds, \quad t>0,
$$
is a time homogeneous diffusion. Other properties of the functional in
(\ref{t}) are studied in the
papers \cite{matsumotoyor99b}, \cite{donatimatsumotoyor00a},
\cite{donatimatsumotoyor00b}, \cite{donatimatsumotoyor00c},
\cite{donatimatsumotoyor00d}. See also Dufresne \cite{dufresne01} and
Matsumoto and Yor  \cite{matsumotoyor03}.

The paper is organized as follows: In the next section a general
method connecting the distribution of a perpetual integral
functional to a hitting time is presented. Examples of the method
are presented in Section \ref{sec41}. In particular, we discuss
the functionals (and their reflected counterparts) appearing  in
(\ref{d-y}) and (\ref{y}). Further, using the Lamperti
transformation we derive the identity (\ref{c-t}) from (\ref{ps2})
and give a new derivation (for another one, see \cite{salminenyor03b}) 
for the joint Laplace transform of the functionals
$$
\int_0^{\infty}
\frac{{\bf
1}_{\{B^{(1/2)}_s> 0\}}}
{(a+\exp(B^{(1/2)}_s))^{2}}\, ds,\quad
L^0_\infty(B^{(1/2)}),\quad
\int_0^{\infty}
\frac{{\bf
1}_{\{B^{(1/2)}_s< 0\}}}
{(a+\exp(B^{(1/2)}_s))^{2}}\, ds,
$$
where $L^0_\infty(B^{(1/2)})$ is the ultimate value of local time
at 0 of $B^{(1/2)}.$ Especially, this latter computation has
interesting connections to some earlier works. We also prove 
in Section \ref{sec41} the identities 
(\ref{i}), (\ref{legall}) and (\ref{hariya}) (see the table below). 
In Section \ref{sec5},
we modify the Feynman-Kac formula to be directly applicable for
computing the Laplace transform of a perpetual functional and
discuss a characterization due to Biane \cite{biane85} for one-sided
functionals. We finish with a short Appendix containing
the Ray--Knight theorems used in the paper.

To summarize the discussion made throughout this work, we use
systematically random time changes in a set-up which a priori
encompasses the scale and speed Feller type representations of
one-dimensional diffusions (for which, see, e.g. the recent paper by
McKean \cite{mckean01}) as well as Lamperti's transformation and
Ray--Knight theorems. Of
course, this stochastic approach and the results it allows to derive agree
with the more analytic Feynman-Kac approach of solving ODE's, and
performing for them the corresponding (deterministic) changes of
variables.

We conclude this introduction with a table containing
(most of) the functionals and the associated hitting
times discussed in this paper. We use the notation
$B^{(\mu)}$ for a Brownian motion with drift $\mu>0,$
$R^{(\delta)}$ for a
Bessel process of dimension $\delta,$ and
$\tilde B^{(\mu)}$ for a reflected Brownian motion with drift $\mu>0.$
All the processes in the table, if nothing else is said,
are started from~0.

$$
\vbox{\offinterlineskip\halign{&\vrule#&\strut\quad\hfil
$#$\quad\hfil\cr
\noalign{\hrule}
height10pt&\omit&&\omit&&\omit&\cr
&{\rm Ref.}&&{\rm Functional}\ (a>0,\, \mu>0)&&{\rm Hitting/occupation\ time}& \cr
height10pt&\omit&&\omit&&\omit&\cr
\noalign{\hrule height0.1pt}
height10pt&\omit&&\omit&&\omit&\cr
&
(\ref{d-y})
&&
\int_0^\infty \exp\bigl(-2aB^{(\mu)}_s\bigr)\,ds
&&
H_0(R^{(2-2\mu/a)}),
&\cr
height10pt&\omit&&\omit&&\omit&\cr
&
&&
&&
\qquad R^{(2-2\mu/a)}_0=1/a
&\cr
height10pt&\omit&&\omit&&\omit&\cr
\noalign{\hrule height0.1pt}
height10pt&\omit&&\omit&&\omit&\cr
&
(\ref{ps2})
&&
\int_0^{\infty}
\exp(-2a\,B^{(\mu)}_s)\,{\bf 1}_{\{B^{(\mu)}_s> 0\}}\, ds
&&
H_{1/a}(R^{(2\mu/a)})
&\cr
height10pt&\omit&&\omit&&\omit&\cr
\noalign{\hrule height0.1pt}
height10pt&\omit&&\omit&&\omit&\cr
&
(\ref{dref2})
&&
\int_0^{\infty}
\exp(-2a\,B^{(\mu)}_s)\,{\bf 1}_{\{B^{(\mu)}_s< 0\}}\, ds
&&
\int_0^\infty {\bf 1}_{\{R^{(2-2\mu/a)}_s>1/a\}}\, ds,
&\cr
height10pt&\omit&&\omit&&\omit&\cr
\noalign{\hrule height0.1pt}
height10pt&\omit&&\omit&&\omit&\cr
&
(\ref{dref3})
&&
\int_0^{\infty}
\exp(2a\,B^{(\mu)}_s)\,{\bf 1}_{\{B^{(\mu)}_s< 0\}}\, ds
&&
\int_0^\infty {\bf 1}_{\{R^{(2+2\mu/a)}_s<1/a\}}\, ds,
&\cr
height10pt&\omit&&\omit&&\omit&\cr
&
&&
&&
\qquad R^{(2+2\mu/a)}_0=1/a
&\cr
height10pt&\omit&&\omit&&\omit&\cr
\noalign{\hrule height0.1pt}
height10pt&\omit&&\omit&&\omit&\cr
&
(\ref{y})
&&
\int_0^{\infty}  (a+\exp(B^{(1/2)}_s))^{-2}\, ds
&&
H_{r}(B^{(a/2)}),\hskip3cm
&\cr
height10pt&\omit&&\omit&&\omit&\cr
&
&&
&&
\qquad r=\frac{1}{a}\log(1+a)%
&\cr
height10pt&\omit&&\omit&&\omit&\cr
\noalign{\hrule height0.1pt}
height10pt&\omit&&\omit&&\omit&\cr
&
(\ref{i})
&&
\int_0^\infty {\bf 1}_{\{B^{(\mu)}_s< 0\}}\, ds
&&
H_\lambda(B^{(\mu)}),\quad \lambda\sim{\rm Exp}(2\mu)
&\cr
height10pt&\omit&&\omit&&\omit&\cr
\noalign{\hrule height0.1pt}
height10pt&\omit&&\omit&&\omit&\cr
&
(\ref{c-t})
&&
\int_0^\infty {\bf 1}_{\{R^{(\delta+2)}_s< 1\}}\, ds
&&
H_1(R^{(\delta)}),\ \delta>0
&\cr
height10pt&\omit&&\omit&&\omit&\cr
\noalign{\hrule height0.1pt}
height10pt&\omit&&\omit&&\omit&\cr
&
(\ref{dm-y})
&&
\int_0^\infty (\exp(2R^{(3)}_s)-1)^{-1}\, ds
&&
H_{\pi/2}(R^{(3)})
&\cr
height10pt&\omit&&\omit&&\omit&\cr
\noalign{\hrule height0.1pt}
height10pt&\omit&&\omit&&\omit&\cr
&
(\ref{legall})
&&
\int_0^\infty \exp(-2R^{(3)}_s)\, ds
&&
H_{1}(R^{(2)})
&\cr
height10pt&\omit&&\omit&&\omit&\cr
\noalign{\hrule height0.1pt}
height10pt&\omit&&\omit&&\omit&\cr
&
(\ref{hariya})
&&
\int_0^\infty (1+R^{(3)}_s)^{-2\gamma}\, ds,
\quad \gamma>1
&&
H_{1/(\gamma-1)}(R^{((2\gamma-1)/(\gamma-1))})
&\cr
height10pt&\omit&&\omit&&\omit&\cr
\noalign{\hrule}}}
$$

\section{Perpetual integral functionals as
first hitting times}
\label{sec1}
Let $\{B^{(\mu)}_t:\ t\geq 0\}$
be a Brownian motion with drift $\mu>0$ and $f$ a non-negative
locally integrable function. Then (see Engelbert and Senf
\cite{engelbertsenf91} and Salminen and Yor \cite{salminenyor03})
\begin{equation}
\label{K1}
I_\infty(f):=\int_0^\infty f(B^{(\mu)}_u) \,  du<\infty \quad {\rm
a.s.}\quad \Leftrightarrow  \int^\infty f(x)\, dx<\infty.
\end{equation}

In our first proposition it is stated, under some additional
assumptions on $f,$ that there exists a diffusion constructed in
the same probability space as $B^{(\mu)}$ such that the perpetual
integral functional $I_\infty(f)$  is a.s. equal to the first hitting
time of a point for this diffusion.



\begin{proposition}
\label{yorprop1} {\sl Let $f:\R\mapsto\R$ be a monotone, twice
continuously differentiable function such that
$r:=\lim_{x\to\infty} f(x)$ exists. Introduce the additive
functional
$$
I_s:=\int_0^s (f'(B^{(\mu)}_u))^2 \,  du.
$$
Assume further that $f'(x)\not= 0$ and
$$
\int^\infty (f'(x))^2 dx <\infty.
$$
Let $Z$ be a diffusion given by
$$
Z_t:=f(B^{(\mu)}_{\alpha_t}),
$$
for $t\geq 0$ such that 
$$
\alpha_t:=\inf\{ s:\, I_s>t\}<\infty.
$$
Then a.s.
\begin{equation}
\label{basic} \int_0^\infty (f'(B^{(\mu)}_s))^2 \,  ds =
\inf\{t:\ Z_t=r\}.
\end{equation}
Moreover,
$Z$ is a solution of the SDE
$$
dZ_t=d\beta_t +
G\circ f^{-1}(Z_t)\, dt,\quad Z_0=f(0),
$$
where $\beta$ is a Brownian motion and 
\begin{equation}
\label{G}
G(x):= {1\over{ (f'(x))^2}}\Bigl({1\over 2}f''(x)+\mu \, f'(x)\Bigr).
\end{equation}
  }
\end{proposition}


\begin{proof}
To fix ideas, take $f$ to be increasing. By Ito's formula
\begin{eqnarray*}
f(B^{(\mu)}_u)-f(0)&=&\int_0^u f'(B^{(\mu)}_s)\, dB^{(\mu)}_s +
{1\over 2} \int_0^u f''(B^{(\mu)}_s)\, ds\\
&=& \int_0^u f'(B^{(\mu)}_s)\, dB_s +
\int_0^u (f'(B^{(\mu)}_s))^2 \, G(B^{(\mu)}_s)\, ds.\\
\end{eqnarray*}
Replacing $u$ by $\alpha_t$ we obtain
$$
Z_t-Z_0=
\int_0^{\alpha_t} f'(B^{(\mu)}_s)\, dB_s +
\int_0^{\alpha_t}  (f'(B^{(\mu)}_s))^2 \, G(B^{(\mu)}_{s})\, ds.
$$
Because
$$
I'_s=(f'(B^{(\mu)}_s))^2
\quad
{\rm and}
\quad
\alpha'_t={1\over {I'_{\alpha_t}}}=(f'(B^{(\mu)}_{\alpha_t}))^{-2}
$$
it follows from L\'evy's theorem that
$$
\beta_t:=
\int_0^{\alpha_t} f'(B^{(\mu)}_s)\, dB_s,
\ t\geq 0,
$$
is a Brownian motion, and we obtain the claimed SDE:
\begin{eqnarray*}
Z_t-Z_0&=&\beta_t +
\int_0^{t}
 (f'(B^{(\mu)}_{\alpha_s}))^2 \,
G(B^{(\mu)}_{\alpha_s})\, d{\alpha_s}\\
&=&\beta_t +
\int_0^{t}G\circ f^{-1}(Z_s)\, ds.
\end{eqnarray*}
Because $B^{(\mu)}_{t}\to\infty$ as $t\to\infty$
we have
$$
\lim_{t\to\infty}f(B^{(\mu)}_{t})=r, \quad {\rm a.s.}
$$
It follows now from
$$
Z_{I_t}
=
f(B^{(\mu)}_{t})< r
$$
by letting $t\to\infty$ that
$$
I_\infty  \
=\ H_r(Z)\quad{\rm a.s.},
$$
completing the proof.
\end{proof}

\begin{remark}
\label{stoppedfunctional}
Let $f$ and $Z$ be as in Proposition \ref{yorprop1}. Then, modifying
the proof above, it is seen for $x>0$ that
\begin{equation}
\label{basic2}
\int_0^{H_x(B^{(\mu)})}(f'(B^{(\mu)}_s))^2 \,  ds =
\inf\{t:\ Z_t=f(x)\}\quad {\rm a.s.}
\end{equation}
\end{remark}

We consider next perpetual integral functionals
restricted on $\R_+.$
For simplicity we take $f$ to be decreasing, and
leave the case ``$f$ increasing'' to the reader.
Further, we use the notation
$$
M^{(-\mu)}_t:=\sup\{B^{(-\mu)}_s:\ s\leq t\} := \sup\{B_s-\mu\,s:\
s\leq t\},\ \mu>0.
$$
Recall that
\begin{equation}
\label{def-rho}
\{\rho^{(\mu)}_t:=M^{(-\mu)}_t-B^{(-\mu)}_t:\ t\geq 0\}
\end{equation}
is identical in law to a reflecting Brownian motion with drift $\mu,$
i.e., the solution of Skorokhod's reflection equation driven by
$$
-B^{(-\mu)}_t=(-B_t)+\mu\,t,\quad t\geq 0.
$$

\begin{proposition} \label{yorprop2} {\sl Let $f$ and $G$ be as in
Proposition \ref{yorprop1}. Assume moreover that $f$ is decreasing and
let $\rho^{(\mu)}$ be as given in (\ref{def-rho}). Let 
$$
A^+_t:=\int_0^t (f'(\rho^{(\mu)}_s))^2 \,  ds,
\quad
{\rm and}
\quad
\alpha^+_t:=\inf\{ s:\, A^+_s>t\}.
$$
Define a process $Z$ via
$$
Z_t:=f(0)-f(\rho^{(\mu)}_{\alpha^+_t})
$$
for $t\geq 0$ such that $\alpha^+_t<\infty.$ Then
\begin{equation}
\label{claim1} I^+_\infty:= \int_0^\infty (f'(B^{(\mu)}_s))^2 \,
{\bf 1}_{\{B^{(\mu)}_s> 0\}}\, ds\ {\mathop=^{\rm{(d)}}}\
\inf\{t:\ Z_t=r^\star\},
\end{equation}
where $r^\star:=f(0)-r.$ The process $Z$ is a solution of the
reflected SDE
$$
dZ_t=d\beta_t -
G\circ f^{-1}(f(0)-Z_t)\, dt + L^0_t(Z) ,\quad Z_0=0,
$$
where $\beta$ is a Brownian motion, 
$\{L^0_t(Z):\ t\geq 0\}$ is a non-decreasing process which increases
only on the zero set of $Z,$ and $G$ is as in Proposition \ref{yorprop1}.
  }
\end{proposition}
\begin{remark}
In comparison with the statement (\ref{basic}) in
Proposition \ref{yorprop1}, notice that the equality in (\ref{claim1})
holds 
``only'' in distribution.
\end{remark}
\begin{proof}
The starting point is the observation 
$$
I^+_\infty  \ {\mathop=^{\rm{(d)}}}\ A^+_\infty,
$$
and, hence, we find the law of $A^+_\infty.$ As in the proof of
Proposition \ref{yorprop1} we use It\^o's formula
\begin{eqnarray*}
&&f(0)-f(\rho^{(\mu)}_t)\\
&&\hskip1cm
=
-\int_0^{t}
f'(\rho^{(\mu)}_s)\,
(dM^{(-\mu)}_s-dB^{(-\mu)}_s)
-\,\frac{1}{2}\,\int_0^{t}
f''(\rho^{(\mu)}_s)\,ds\\
&&
\hskip1cm
=
\int_0^{t}
f'(\rho^{(\mu)}_s)\, dB_s
-f'(0) M^{(-\mu)}_t
-\int_0^{t} \Big(\frac{1}{2}f''(\rho^{(\mu)}_s)
 +\mu f'(\rho^{(\mu)}_s)\Big)\, ds.
\end{eqnarray*}
Consequently, for $t\geq 0$ such that $\alpha^+_t<\infty$
\begin{eqnarray*}
&&Z_t
=
\int_0^{\alpha^+_t}
f'(\rho^{(\mu)}_s)\, dB_s
-f'(0) M^{(-\mu)}_{\alpha^+_t}
-\int_0^{\alpha^+_t} \big(f'(\rho^{(\mu)}_s)\big)^2
\ G(\rho^{(\mu)}_s)\, ds\\
&&
\hskip.5cm
= \beta_t-
\int_0^{t}G\circ f^{-1}(f(0)-Z_s)\, ds+L^0_t(Z),
\end{eqnarray*}
where $L^0_t(Z):=-f'(0) M^{(-\mu)}_{\alpha^+_t}$ determines a
non-decreasing process which increases only on the zero set of
$\{Z_t:\ t\geq 0\}.$ It holds also that
$$
0\leq Z_t< r^*:=f(0)-r
$$
for all $t\geq 0$ such that $\alpha^+_t<\infty.$ To conclude the
proof, observe that
$$
Z_{A^+_\infty}=\lim_{t\to\infty}\Big(
f(0)-f(\rho^{(\mu)}_t)\Big)=r^*
$$
giving (\ref{claim1}).
\end{proof}
\begin{remark}
The method of random time change is also used by
Getoor and Sharpe \cite{getoorsharpe79}
Section 5  as a tool to
compute the distributions of stopped functionals of Brownian
motion. Their idea is similar to the one in Propositions
\ref{yorprop1} and \ref{yorprop2}, that is, to identify the
functional with the first hitting time of a point for some other
diffusion.
\end{remark}

\section{Detailed studies of some perpetual functionals}%
\label{sec41}

\subsection{Dufresne's functional}
\label{ct-revis}
As a very particular case of our method, we study now the identity
(\ref{d-y}) (or (\ref{d-y2}) below)
and its various one-sided variants. For more
details concerning (\ref{d-y}),
we refer to \cite{yor01} p. 16. For the joint distribution of
the functionals in (\ref{332}) and (\ref{333}), see Salminen and Yor
\cite{salminenyor03b}. In Section 4 Example 4.7 an additional
characterization is derived for the functionals in (\ref{333}).
\begin{proposition}
\label{refdufresne}
{\sl Let $B^{(\mu)}$ with $\mu>0$ be started from
0. Then for $a>0$ the following 5 identities hold}
\hfill\break\hfill
{\bf a:}\ (cf. (\ref{d-y})),
{\sl
\begin{equation}
\label{d-y2}
\int_0^\infty \exp(-2aB_s^{(\mu)})\, ds
\quad{\mathop=^{\rm{(d)}}}\quad H_0(R^{(2-2\mu/a)}),
\end{equation}
where $R^{(2-2\mu/a)}$ is a Bessel process
of dimension $2-2\mu/a$
started at 1/a.}
\hfill\break\hfill
{\bf b:}\ (cf. (\ref{ps2})),
{\sl
\begin{equation}
\label{331}
\int_0^{\infty} \exp(-2a\,B^{(\mu)}_s)\,
 {\bf 1}_{\{B_s^{(\mu)}> 0\}}
\, ds
\quad{\mathop=^{\rm{(d)}}}\quad
H_{1/a}(R^{\,(2\mu/a)}),
\end{equation}
where $R^{\,(2\mu/a)}$ is started at 0 and, in
the case $0<\mu<a,$ reflected at 0.
}
\hfill\break\hfill
{\bf c:}\ (cf. (\ref{ps2})),
{\sl
\begin{equation}
\label{332}
\int_0^{\infty}
\exp(-2a\,B^{(\mu)}_s)\,
{\bf 1}_{\{B^{(\mu)}_s> 0\}}\, ds
\quad{\mathop=^{\rm{(d)}}}\quad
\int_0^\infty
{\bf 1}_{\{R^{(2-2\mu/a)}_{s}< 1/a\}}\, ds,
\end{equation}
where $R^{(2-2\mu/a)}$ is started at $1/a$ and
killed when it hits 0.
}
\hfill\break\hfill
{\bf d:}\ (cf. (\ref{dref2}) and (\ref{b11})),
{\sl
 \begin{equation}
\label{333}
\int_0^{\infty}
\exp(-2a\,B^{(\mu)}_s)\,
\,
{\bf 1}_{\{B^{(\mu)}_s< 0\}}\, ds
\quad{\mathop=^{\rm{(d)}}}\quad
\int_0^\infty
{\bf 1}_{\{R^{(2-2\mu/a)}_{s}> 1/a\}}\, ds,
\end{equation}
where $R^{(2-2\mu/a)}$ is started from $1/a$ and
killed when it hits 0.
}
\hfill\break\hfill
{\bf e:}\ (cf. (\ref{dref3})),
{\sl
 \begin{equation}
\label{334}
\int_0^{\infty}
\exp(2a\,B^{(\mu)}_s)\,
\,
{\bf 1}_{\{B^{(\mu)}_s<0\}}\, ds
\quad{\mathop=^{\rm{(d)}}}\quad
\int_0^{\infty}
{\bf 1}_{\{R^{(2+2\mu/a)}_s<1/a\}}\, ds,
\end{equation}
where $R^{(2+2\mu/a)}_0=1/a.$
}

\end{proposition}

\begin{remark}
In fact, as is clear from the proofs below, the identities
(\ref{d-y2}), (\ref{332}), (\ref{333}), and (\ref{334}) hold a.s.,
i.e., the Bessel processes appearing therein can be
constructed in the same probability space where $B^{(\mu)}$ is given.
\end{remark}
\begin{proof}
\hfill\break\hfill
{\bf a:}\quad We let $a=1;$ the case with arbitrary $a$ can be treated
using the scaling property (however, see the proof of (\ref{331})
where we work with arbitrary $a$). Apply Proposition \ref{yorprop1} with
$f(x)=\exp(-x).$ Now
$f$ is decreasing, and straightforward computations show that
$$
G(x)=(\frac{1}{2}-\mu)\exp(x).
$$
The process $Z$ is given by the SDE
$$
d\, Z_t=d\, \beta_t+\frac{1-2\mu}{2\,Z_t}\, dt,\quad Z_0=f(0)=1.
$$
Consequently, $Z$ is a Bessel process of dimension
$\delta=2-2\mu,$ as claimed in (\ref{d-y2}).

\hfill\break\hfill
{\bf b:}\quad Here we use Proposition \ref{yorprop2} with
$f(x)=a^{-1}\, \exp(-ax)$  which leads us to
consider the additive functional
$$
A^+_t:= \int_0^{t} \exp(-2a\,\rho^{(\mu)}_s)\, ds
$$
where  $\rho^{(\mu)}$ is a reflecting Brownian motion with drift $\mu>0$ (see
(\ref{def-rho})). By Proposition \ref{yorprop2} the process $Z$
given by
$$
Z_t:=\frac{1}{a}\big( 1-\exp(-a\,\rho^{(\mu)}_{\alpha^+_t}) \big)
$$
satisfies the SDE
$$
Z_t=\beta_t+ (\frac{\mu}{a}-\frac{1}{2})\int_0^{t}
\big(\frac{1}{a}-Z_s\big)^{-1}\, ds +L^0_t(Z),\ Z_0=0.
$$
Recall that $\alpha^+$ is the inverse of $A^+,$ and $L^0(Z)$
is a non-decreasing process which increases only on the
zero set of $Z.$ Proposition \ref{yorprop2} yields now
\begin{eqnarray*}
&&\int_0^{\infty}
\exp(-2a\,B^{(\mu)}_s)\,
 {\bf 1}_{\{B_s^{(\mu)}> 0\}}\, ds
\quad{\mathop=^{\rm{(d)}}}\quad
\int_0^{\infty}
\exp(-2a\,\rho^{(\mu)}_s)\, ds\\
&&\hskip5.5cm
\quad{\mathop=^{\rm{(d)}}}\quad
H_{1/a}(Z).
\end{eqnarray*}
It remains to show (cf. (\ref{331})) that  $H_{1/a}(Z)$ is identical in law with
$H_{1/a}(R^{(2\mu/a)}).$ One way to do this is to compute
the Laplace transform of $H_{1/a}(Z)$ by using the standard
diffusion theory (see It\^o and McKean
\cite{itomckean74} or \cite{borodinsalminen02} II.10 p. 18). A
slightly shorter
way is to observe that $\{ 1/a -Z_t:\ t\geq 0\}$ is a Bessel process
of dimension $2-2\mu/a$ started at $1/a,$ reflected at $1/a$ and
killed when it hits 0, denoted $\tilde R^{(2-2\mu/a)}$.
This gives us
$$
H_{1/a}(Z)
\quad{\mathop=^{\rm{(d)}}}\quad
H_0(\tilde R^{(2-2\mu/a)}).
$$
Finally, the claim
\begin{equation}
\label{eq21}
H_0(\tilde R^{(2-2\mu/a)})
\quad{\mathop=^{\rm{(d)}}}\quad
H_{1/a}(R^{(2\mu/a)}),
\end{equation}
where the Bessel process $R^{(2\mu/a)}$ is started at $0,$ can be
verified by straightforward but lenghty computations with Laplace
transforms. (As seen in Section 3.2, the Ciesielski--Taylor
identity follows easily from (\ref{eq21})).

\hfill\break\hfill
{\bf c} \& {\bf d:}\quad To prove (\ref{332}) and (\ref{333}) we use the
Lamperti transformation
(see Lamperti \cite{lamperti72}) as in the
proof of Dufresne's identity in Yor
\cite{yor92b}. Notice that the Lamperti transformation is used
(and proved) implicitly also in the above proof of (\ref{d-y2}). To recall the
transformation, let
$$
A^{(1)}_t=\int_0^{t}
\exp(-2a\,B^{(\mu)}_s)\,ds
$$
and
$$
\alpha^{(1)}_t:=
\inf\{s:\ A^{(1)}_s>t\}.
$$
Then the (first) Lamperti transformation says that the process
\begin{equation}
\label{lam1}
\{ \frac{1}{a}\exp(-a\,B^{(\mu)}_{\alpha^{(1)}_t}):\ t\geq 0\}
\end{equation}
is a Bessel process of dimension $2-2\mu/a$
started at $1/a$ and
killed when it hits 0. Letting $R^{(2-2\mu/a)}$
denote this Bessel process we obtain
\begin{eqnarray}
\label{psalt1} && \nonumber \int_0^{\infty}
\exp(-2a\,B^{(\mu)}_s)\,
{\bf 1}_{\{B^{(\mu)}_s> 0\}}\, ds\\
&&\hskip3cm\nonumber
=\int_0^{\infty}
\exp(-2a\,B^{(\mu)}_{\alpha^{(1)}_s})\,
{\bf 1}_{\{B^{(\mu)}_{\alpha^{(1)}_s}\,> 0\}}
\,
d\alpha^{(1)}_s\\
&&
\hskip3cm
=
\int_0^{H_0(R^{(2-2\mu/a)})}
{\bf 1}_{\{R^{(2-2\mu/a)}_{s}< 1/a\}}\, ds\\
\nonumber 
&&\hskip3cm
 =
\int_0^\infty
{\bf 1}_{\{R^{(2-2\mu/a)}_{s}< 1/a\}}\, ds.
\end{eqnarray}
For (\ref{333}) we have similarly
\begin{eqnarray}
\label{psalt2} && \nonumber\int_0^{\infty}
\exp(-2a\,B^{(\mu)}_s)\,
{\bf 1}_{\{B^{(\mu)}_s< 0\}}\, ds\\
&&\hskip3cm\nonumber
=
\int_0^{\infty}
\exp(-2a\,B^{(\mu)}_{\alpha^{(1)}_s})\,
{\bf 1}_{\{B^{(\mu)}_{\alpha^{(1)}_s}< 0\}}
\, d\alpha^{(1)}_s
\\
&&
\hskip3cm
=\int_0^{H_0(R^{(2-2\mu/a)})} {\bf 1}_{\{R^{(2-2\mu/a)}_{s}>
1/a\}}\, ds.
\\
\nonumber
&&\hskip3cm =
\int_0^\infty
{\bf 1}_{\{R^{(2-2\mu/a)}_{s}> 1/a\}}\, ds.
\end{eqnarray}

\hfill\break\hfill
{\bf e:}\quad  The identity (\ref{334}) can be proved using the
Lamperti transformation (\ref{lam1}) with negative $\mu.$
This can also be formulated for $\mu>0$ by defining
$$
A^{(2)}_t:=\int_0^{t}
\exp(2a\,B^{(\mu)}_s)\,ds
$$
and
$$
\alpha^{(2)}_t:=
\inf\{s:\ A^{(2)}_s>t\}.
$$
Then the (second) Lamperti transformation states that
\begin{equation}
\label{lam2}
\{ \frac{1}{a}\exp(a\,B^{(\mu)}_{\alpha^{(2)}_t}):\ t\geq 0\}
\end{equation}
is a Bessel process of dimension $2+2\mu/a$, denoted
$R^{(2+2\mu/a)},$ started
at $1/a.$ Consequently, as in the proofs of {\bf c} and {\bf d} above, 
\begin{eqnarray*}
&&\int_0^{\infty}
\exp(2a\,B^{(\mu)}_s)\,
\,
{\bf 1}_{\{B^{(\mu)}_s<0\}}\, ds
=
\int_0^{\infty}
{\bf 1}_{\{R^{(2+2\mu/a)}_s<1/a\}}\, ds.
\end{eqnarray*}
\end{proof}
\begin{remark}
\label{r1}

{\bf (1)} Notice that in case $a=1$ and $\mu=1/2$ the identity
(\ref{331}) (see also (\ref{psalt1})) takes the form
\begin{eqnarray*}
&&\int_0^{\infty} \exp(-2\,B^{(1/2)}_s)\,
 {\bf 1}_{\{B_s^{(1/2)}> 0\}}
\, ds \quad{\mathop=^{\rm{(d)}}}\quad H_{1}(\rho)\\
&&\hskip5.7cm \quad{\mathop=^{\rm{(d)}}}\quad \int_0^{H_0(B)} {\bf
1}_{\{B_s<1\}} \, ds
\end{eqnarray*}
 where $\rho$ is a reflecting Brownian motion started at 0. Similarly,
(\ref{psalt2}) can be written as
\begin{eqnarray*}
&&\int_0^{\infty} \exp(-2\,B^{(1/2)}_s)\,
 {\bf 1}_{\{B_s^{(1/2)}< 0\}}
\, ds
=
\int_0^{H_0(B)} {\bf 1}_{\{B_s>1\}} \, ds.
\end{eqnarray*}
From \cite{borodinsalminen02} formula 1.2.4.1 p. 200, for $\gamma\geq
0$
$$
\E_1\Big( \exp\Big(-\gamma \int_0^{H_0(B)} {\bf 1}_{\{B_s>1\}} \,
ds\Big)\Big) =\frac{1}{1+\sqrt{2\gamma}}.
$$
By the well known formula (cf. \cite{borodinsalminen02} formula
1.2.0.1 p. 198)
$$
\E_x\Big(\exp(-\gamma H_0(B))\Big)=\exp(-x\sqrt{2\gamma}),\quad x>0,\
\gamma\geq 0,
$$
and, therefore,
\begin{eqnarray*}
&&\int_0^{\infty} \exp(-2\,B^{(1/2)}_s)\,
 {\bf 1}_{\{B_s^{(1/2)}< 0\}}
\, ds
 \quad{\mathop=^{\rm{(d)}}}\quad
H_\lambda(B),
\end{eqnarray*}
where $B$ is started at 0 and $\lambda$ is an exponentially (with
parameter 1) distributed random variable independent of $B.$

\hfill\break\hfill
{\bf (2)}\hskip4mm By the Lamperti time
changes (\ref{lam1}) and (\ref{lam2}), a general 
perpetual integral
functional of geometric Brownian motion can be expressed in terms of
Bessel processes as follows:
\begin{eqnarray*}
&&\int_0^{\infty}
f(\exp(a\,B^{(\mu)}_s))\, ds
=
\int_0^{\infty}(a\,R^{(2+2\mu/a)}_s)^{-2}\ f(a\,R^{(2+2\mu/a)}_s)\ ds\\
&&\hskip4cm
=
\int_0^{\infty}(a\,R^{(2-2\mu/a)}_s)^{-2}\
f\Big(\,{1\over{a\,R^{(2-2\mu/a)}_s}}\,\Big)\ ds,
\end{eqnarray*}
where $R^{(2-2\mu/a)}$ is killed when it hits 0.
\end{remark}

\subsection{Ciesielski--Taylor identity}

The identity (\ref{331}) (cf. also (\ref{eq21})) is now applied to deduce
the well known and puzzling Ciesielski--Taylor identity
(\ref{c-t}) (or (\ref{c-t1}) below), see Ciesielski and Taylor
\cite{ciesielskitaylor62}, and  also Yor \cite{yor91} and \cite{yor92c}, Chap. IV.
\begin{proposition}
{\sl The following identity in law holds:
\begin{equation}
\label{c-t1}
\int_0^\infty {\bf 1}_{\{R^{(\delta+2)}_s< 1/a\}}\, ds
\quad{\mathop=^{\rm{(d)}}}\quad
H_{1/a}(R^{(\delta)}),\quad
\forall\ \delta>0
\end{equation}
where the Bessel processes are started at 0.
}
\end{proposition}
\begin{proof} By Proposition \ref{refdufresne} {\bf b} and {\bf c}
$$
H_{1/a}(R^{(2\mu/a)})
\quad{\mathop=^{\rm{(d)}}}\quad
\int_0^\infty {\bf 1}_{\{R^{(2-2\mu/a)}_t<1/a\}}\ dt,
$$
where $R^{(2-2\mu/a)}$ is started from $1/a.$ Using the standard time
reversal  argument (see Pitman and Yor  \cite{pitmanyor81} p. 341)
we obtain
$$
\int_0^\infty {\bf 1}_{\{R^{(2-2\mu/a)}_t<1/a\}}\ dt
\quad{\mathop=^{\rm{(d)}}}\quad
\int_0^\infty {\bf 1}_{\{R^{(2+2\mu/a)}_t<1/a\}}\ dt,
$$
where $R^{(2+2\mu/a)}$ is started at $0.$ Letting $\delta=2\mu/a$
gives now the claimed identity (\ref{c-t1}).
\end{proof}

\subsection{Translated Dufresne's functional}

An interesting case emerging from Proposition \ref{yorprop1} is
when the function $G$ in (\ref{G}) is equal to a constant $c,$
say. Straightforward computations show that this yields us the
functional
$$
\int_0^\infty \bigl({c\over\mu}+A \exp(2\mu B^{(\mu)}_s)\bigr)^{-2}\, ds,
$$
where $A$ is a free constant. It is natural to assume
that $c>0$ and $A>0.$ By the scaling
property of Brownian motion,
$$
\int_0^\infty \bigl({c\over\mu}+A \exp(2\mu B^{(\mu)}_s)\bigr)^{-2}\,
ds
\ {\mathop=^{\rm{(d)}}}\
{1\over{4\mu^2}}\int_0^\infty \bigl({c\over\mu}+
A \exp(B^{(1/2)}_s)\bigr)^{-2}\, ds
$$
indicating that we are, in fact, obtaining information only about
Brownian motion with drift $\mu=1/2.$ Choosing the constants appropriately,
this functional, which in \cite{salminenyor03b} is called
Dufresne's translated perpetuity, can be expressed in the form
$$
\int_0^{\infty}  (a\,\exp(B^{(1/2)}_s)+1)^{-2}\, ds
$$
where $a>0.$ In \cite{salminenyor03b} we derive the joint Laplace
transform for the one sided variants of this functional, that is,
for the functionals in (\ref{342}) and (\ref{343}) below. Here we give a
new derivation of this Laplace transform based on the Ray--Knight
theorem for Brownian motion stopped at the first hitting time. However, before
this, we characterize the functional and its one-sided variants
via hitting times.

\begin{proposition}
\label{trandufresne} {\sl Suppose that $B^{(1/2)}_0=0.$ Then the
following 3 identities hold} \hfill\break\hfill {\bf a:}
{\sl
\begin{equation}
\label{341} \int_0^{\infty}  (a\exp(B^{(1/2)}_s)+1)^{-2}\, ds
\quad{\mathop=^{\rm{(d)}}}\quad H_{r}(\beta^{(1/2)}),
\end{equation}
where $a>0,$ $r=\log((1+a)/a),$ and $\beta^{(1/2)}$ is a Brownian
motion with drift $1/2$ started at 0.
}
\hfill\break\hfill {\bf b:} {\sl
\begin{equation}
\label{342} \int_0^\infty \Big(a\exp(B^{(1/2)}_s)+1\Big)^{-2} {\bf
1}_{\{B^{(1/2)}_s>0\}} ds \quad{\mathop=^{\rm{(d)}}}\quad
H_{r}(\tilde \beta^{(1/2)}),
\end{equation}
where $a$ and $r$ are as above, and  $\tilde \beta^{(1/2)}$ is
reflecting Brownian motion with drift $1/2$ started at 0. 
} \hfill\break\hfill {\bf c:} {\sl
\begin{equation}
\label{343}
\int_0^\infty \Big(a\exp(B^{(1/2)}_s)+1\Big)^{-2}
{\bf 1}_{\{B^{(1/2)}_s<0\}} ds
\ {\mathop=^{\rm{(d)}}}\
H_\lambda(\beta^{(1/2)}),
\end{equation}
where $a>-1,$ $\lambda$ is exponentially
distributed with parameter $(1+a),$ and $\beta^{(1/2)}$ is as in (\ref{341}).
}
\end{proposition}
\begin{proof}
\hfill\break\hfill {\bf a:}\quad
Instead of simply refering to Proposition \ref{yorprop1},
it is perhaps more instructive to go through the computation; in
this way we also gain better understanding why the case $\mu=1/2$
is a special one. Define
$$
f(x):= \log\Big({a\over{a+\exp(-x)}}\Big),
$$
and observe that
$$
f'(x)= (a\exp(x)+1)^{-1}.
$$
By It\^o's formula
%
\begin{eqnarray}
\label{e00}
\nonumber
&&\log\big(a+\exp(-B^{(\mu)}_t)\big)-\log(a+1)
=-\int_0^{t}  {d B_s\over{a\exp(B^{(\mu)}_s)+1}}\hskip3cm\\
&&\hskip1cm
+({1\over 2}-\mu)\int_0^{t}  {d s\over{a\exp(B^{(\mu)}_s)+1}}
-{1\over 2}\int_0^{t}  {d s\over{(a\exp(B^{(\mu)}_s)+1)^{2}}}.
\end{eqnarray}
Introduce
$$
A_t:=\int_0^{t}  {d s\over{(a\exp(B^{(\mu)}_s)+1)^{2}}},
$$
and
$$
Z_u:=f(B^{(\mu)}_{\alpha_u})=\log\Big({a\over{a+
\exp(-B^{(\mu)}_{\alpha_u})}}\Big),
$$
where $\alpha$ is the inverse of $A.$ We have
$Z_0=\log(a/(1+a))<0,$ and
$Z_u<0$ for all $u$ such that $\alpha_u<\infty.$ 
Moreover, from (\ref{e00}) it is seen that $Z$ is a solution of the SDE
$$
Z_u=Z_0+\beta^{(1/2)}_u + (\mu-{1\over 2})\int_0^{u}  {d
s\over{1-\exp(Z_s)}},
$$
where $\beta^{(1/2)}$ is a Brownian motion with drift $1/2$ starting
from 0. Consequently, because $Z_u\to 0$ as $u\to\infty$ we obtain,
as explained in the proof of Proposition
\ref{yorprop1},
$$
\int_0^{\infty}  {d s\over{(a\exp(B^{(\mu)}_s)+1)^{2}}}
=
H_0(Z),
$$
and, in particular, for $\mu=1/2$ we have (\ref{341}).
\hfill\break\hfill {\bf b:}\quad The identity (\ref{342}) can be
proved by directly applying Proposition \ref{yorprop2}, and
we leave the details to the reader. 
\hfill\break\hfill {\bf
c:}\quad Consider now the identity (\ref{343}). Our proof of this
relies on the Ray--Knight Theorem \ref{r-k2} given in Appendix at
the end of the paper. Firstly, for arbitrary $\mu>0$ we have by
the occupation time formula
$$
\int_0^\infty \Big(a\exp(B^{(\mu)}_s)+1\Big)^{-2} {\bf
1}_{\{B^{(\mu)}_s<0\}} ds =\int_{-\infty}^0
\Big(a\exp(y)+1\Big)^{-2}\,L^{y}_\infty(B^{(\mu)})\, dy,
$$
where  $L^{y}_\infty(B^{(\mu)})$ is the total local time of
$B^{(\mu)}$ at level $y.$ By Theorem \ref{r-k2},
$$
\int_{-\infty}^0
\Big(a\exp(y)+1\Big)^{-2}\,L^{y}_\infty(B^{(\mu)})\, dy,
\quad{\mathop=^{\rm{(d)}}}\quad \int_{0}^\infty
\Big(a\exp(-y)+1\Big)^{-2}\,Z_{y}\, dy,
$$
where $Z$ satisfies the SDE
$$
    dZ_y= 2\sqrt{Z_y}\,dB_y -2\mu Z_y\,dy.
$$
The distribution of $Z_0$ is the distribution of
$L^0_\infty(B^{(\mu)})$, i.e., the exponential distribution with
parameter $\mu.$ Consider
\begin{eqnarray}
\label{345} \nonumber
 &&(a\,\exp({-y})+1)^{-1}\,
Z_y-(a+1)^{-1}\,Z_0 \\
\nonumber
&& \hskip2cm =2 \int_0^y (a\,\exp({-u})+1)^{-1}\,\sqrt{Z_u}\,dB_u \\
&&\hskip4cm - \int_0^y
\frac{2\mu+(2\mu-1)\,a\,\exp({-u})}{(a\,\exp({-u})+1)^{-2}}\,Z_u\,du.
\end{eqnarray}
Define
$$
C_y:=\int_0^y (a\,\exp({-u})+1)^{-2}\, Z_u\, du, \quad y\geq 0,
$$
and let $c$ denote the inverse of $C.$ From
(\ref{345}) we obtain for $\mu=1/2$
$$
(a\,\exp({-c_y})+1)^{-1}\, Z_{c_y}-(a+1)^{-1} \,Z_0
=2\,\beta_{y}-{y},
$$
where $\beta$ is a Brownian motion. It follows, because $Z_t\to 0$ as
$t\to\infty,$ that
$$
\int_0^\infty \Big(a\exp(B^{(1/2)}_s)+1\Big)^{-2} {\bf
1}_{\{B^{(1/2)}_s<0\}} ds \ {\mathop=^{\rm{(d)}}}\ \inf\{ y:\
\beta^{(1/2)}_y=\xi\},
$$
where $\xi=Z_0/2(a+1)$ is exponentially distributed with parameter
$(a+1).$
\end{proof}

\begin{remark}
The formula (\ref{basic2}) in Remark 2.2 gives 
$$
\int_0^{H_p(B^{(1/2)})}
{d s\over{(a\exp(B^{(1/2)}_s)+1)^{2}}}
=
H_q(\beta^{(1/2)}),
$$
where $p>0,$ and $q =\log\Big((a+1)/(a+\exp({-p}))\Big).$
\end{remark}

We proceed by computing the joint Laplace transform of
the functionals appearing in (\ref{342}) and (\ref{343}).
This Laplace transform is also given in \cite{salminenyor03b} but here
we give a new derivation which from our point of view
has independent interest. Define
$$
\Delta^{(\pm)}_{a}:=\int_0^{\infty} {{\bf
1}_{\{B_s^{(1/2)}\in\R_\pm\}}\over {(a\exp(B^{(1/2)}_s)+1)^{2}}}\ ds.
$$

\begin{proposition}
\label{T2} {\it For non-negative $k, K$ and $c$
\begin{eqnarray}
\label{F} &&  F(k,c,K):= \E_0\Big(\exp \Big(-k\,
\Delta^{(+)}_{a}-{c}\, L^0_\infty(B^{(1/2)}) -{K}\,
 \Delta^{(-)}_{a}\Big)\Big)\\&&\hskip.5cm
\nonumber={\sqrt{8k+1}\,\exp({r\over 2})\over{ \sqrt{8k+1}\,{\rm
cosh}({r\over 2}\,\sqrt{8k+1}) + (2c(a+1)+\sqrt{8K+1})\,{\rm
sinh}({r\over 2}\,\sqrt{8k+1})}}.
\end{eqnarray}
where $r=\log((a+1)/a).$ In particular,
with the notation as in Proposition \ref{trandufresne},
\begin{eqnarray}
\label{F1}
\nonumber
&&
F(k,0,k)=
\E_0\Big(\exp\Big(-k\,\int_0^{\infty}  (a\exp(B^{(1/2)}_s)+1)^{-2}\,
ds\Big)\Big)\\
\nonumber
&&\hskip1.8cm
=\E_0\Big(\exp\Big(-k\,H_r(\beta^{(1/2)})\Big)\Big)\\
\nonumber
&&\hskip1.8cm
\\
\nonumber
&&\hskip1.9cm
=\exp\Big(-{r\over 2}(\sqrt{8k+1}-1)\Big),
\end{eqnarray}
\begin{eqnarray}
\label{F2}&&
\nonumber
F(0,0,K)=
\E_0\Big(\exp\Big(-K\,\int_0^{\infty} {{\bf 1}_{\{B_s^{(1/2)}<
0\}}\over {(a\exp(B^{(1/2)}_s)+1)^{2}}}\, ds\Big)\Big)
\\
\nonumber
&&\hskip2cm
=\E_0\Big(\exp\Big(-K\,H_\lambda(\beta^{(1/2)})\Big)\Big)
\\
\nonumber
&&\hskip2cm
=\frac{1+a}{a+\frac{1}{2}+\sqrt{2K+\frac{1}{4}}},
\end{eqnarray}
and 
\begin{eqnarray}
\label{I+} &&\nonumber
F(k,0,0)=
\E_0\Big(\exp\Big(-k\,\int_0^{\infty} {{\bf 1}_{\{B_s^{(1/2)}>
0\}}\over {(a\exp(B^{(1/2)}_s)+1)^{2}}}\, ds\Big)\Big)
\\
\nonumber
&&\hskip1.8cm
=\E_0\Big(\exp\Big(-k\,H_r(\tilde\beta^{(1/2)})\Big)\Big)
\\
\nonumber
&&\hskip1.8cm=
{\sqrt{8k+1}\,\exp({r\over 2})\over{
\sqrt{8k+1}\,{\rm cosh}({r\over 2}\,\sqrt{8k+1}) + {\rm
sinh}({r\over 2}\,\sqrt{8k+1})}}.
\end{eqnarray}
Moreover, $\Delta^{(+)}_a$ and $\Delta^{(-)}_a$ are conditionally
independent given $L^0_\infty(B^{(1/2)}).$ 
}
\end{proposition}
\begin{proof}
We begin as in \cite{salminenyor03b} and express
the functionals  $\Delta^{(\pm)}_{a}$ in terms of a Brownian motion
without drift. Firstly, notice that
$$ \Delta^{(+)}_{a}
\quad{\mathop=^{\rm{(d)}}}\quad \int_0^{\infty}
{\exp(2\,B^{(-1/2)}_s)\,{\bf 1}_{\{B_s^{(-1/2)}< 0\}}\over
{(a+\exp(B^{(-1/2)}_s))^{2}}}\ ds,
$$
and
$$
\Delta^{(-)}_{a}
\quad{\mathop=^{\rm{(d)}}}\quad \int_0^{\infty}
{\exp(2\,B^{(-1/2)}_s)\,{\bf 1}_{\{B_s^{(-1/2)}> 0\}}\over
{(a+\exp(B^{(-1/2)}_s))^{2}}}\ ds.
$$
We apply the Lamperti transformation (\ref{lam2}) with $\mu=-1/2$ and
$a=1$ on
the left hand sides of these identities. Recalling that $R^{(1)}$
is, in fact, a Brownian motion killed
when it hits 0 we obtain
$$
\Delta^{(+)}_{a} \quad{\mathop=^{\rm{(d)}}}\quad \int_0^{H_0(B')} {{\bf
1}_{\{B'_s< 1\}}\over {(a+B'_s)^{2}}}\ ds,
\quad
\Delta^{(-)}_{a}
\quad{\mathop=^{\rm{(d)}}}\quad \int_0^{H_0(B')} {{\bf
1}_{\{B'_s> 1\}}\over {(a+B'_s)^{2}}}\ ds,
$$
where $B'$ is a Brownian motion started at 1. 
Instead of working with a Brownian motion $B'$ starting at 1 we
introduce a ``new'' Brownian motion $B=1-B'$ starting at 0. Using
spatial homogeneity of Brownian motion we obtain
$$
\int_0^{H_0(B')}
{{\bf 1}_{\{B'_s< 1\}}\over
{(a+B'_s)^{2}}}\ ds
\quad{\mathop=^{\rm{(d)}}}\quad
\int_0^{H_1(B)}
{{\bf 1}_{\{B_s> 0\}}\over
{(a+1-B_s)^{2}}}\ ds=:\Delta^+,
$$
and
$$
\int_0^{H_0(B')}
{{\bf 1}_{\{B'_s> 1\}}\over
{(a+B'_s)^{2}}}\ ds
\quad{\mathop=^{\rm{(d)}}}\quad
\int_0^{H_1(B)}
{{\bf 1}_{\{B_s< 0\}}\over
{(a+1-B_s)^{2}}}\ ds=:\Delta^-.
$$
By the occupation time formula
\begin{eqnarray}
\label{variables2} \nonumber &&\big(\,\Delta^{(+)}_{a},\
L^0_\infty(B^{(1/2)}),
\ \Delta^{(-)}_{a}\,\big)
\\
\nonumber
&&
\hskip1.5cm
\quad{\mathop=^{\rm{(d)}}}\quad
\big(\,\Delta^+,\ L^0_{H_1}(B),
\ \Delta^-\,\big)
\\
\nonumber
&&
\hskip2.0cm
=
\Big(\,\int_0^1
{L^{1-y}_{H_1}(B)\over
{(a+y)^{2}}}\,dy\ ,\
L^0_{H_1}(B)\ ,\
\int_1^\infty
{L^{1-y}_{H_1}(B)\over
{(a+y)^{2}}}\,dy\,\Big),
\end{eqnarray}
where $L^{1-y}_{H_1}(B)$ is the local time of $B$ at $1-y$ up to
$H_1(B).$ From the proof of (\ref{343}) in Proposition
\ref{trandufresne} we know
\begin{eqnarray}
\label{344}
\nonumber
&&\E_0\Big(\exp\Big(-K\,\Delta^{-}\Big)
\ |\ L^0_{H_1}(B)=u\Big)\\
\nonumber
&&\hskip2cm
=\E_0\Big(\exp\Big(-K\,\Delta^{(-)}_a\Big)
\ |\ L^0_{\infty}(B^{(1/2)})=u\Big)
\\
\nonumber
&&\hskip2cm
=\exp\Big(-\frac{u}{4(1+a)}(\sqrt{8\,K+1}
-1)\Big).
\end{eqnarray}
From the Ray--Knight Theorem \ref{r-k1} it now follows
\begin{eqnarray}
\label{variables2}
\nonumber
&&\E_0\Big(\exp\Big(\,-k\,\Delta^{(+)}_{a}-c\, L^0_\infty(B^{(1/2)})
-K\, \Delta^{(-)}_{a}\,\Big)\Big)\\
\nonumber
&&
\hskip0.5cm
=
\E_0\Big(\exp\Big(\,-k\,
\int_0^1
{X^{(2)}_s\over
{(a+s)^{2}}}\,ds\, -
c\,X^{(2)}_1
-{X^{(2)}_1\over{4\,\alpha}}( \sqrt{8K+1}-1)\Big)\Big)\\
&&
\hskip0.5cm
=
\E_0\Big(\exp\Big(\,-k\,
\int_0^1
{X^{(2)}_s\over
{(a+s)^{2}}}\,ds\, -
\gamma\,X^{(2)}_1\Big)\Big),
\end{eqnarray}
where $X^{(2)}$ denotes the 2-dimensional squared Bessel process
started at 0 and
\begin{equation}
\label{gamma}
\gamma=c+{1\over{4\,(a+1)}}( \sqrt{8K+1}-1).
\end{equation}
Next, recall from Pitman and Yor \cite{pitmanyor82} or Revuz and
Yor \cite{revuzyor01} Exercise 1.34 p. 453 that for a general
squared Bessel process $X^{(\delta)},\ \delta\geq 2,$
and for any positive 
Radon measure $m$
on $[0,\infty),$ we have for $x\geq 0$ and $t>0$
\begin{eqnarray*}
&&
\E_x\Big(\exp\Big(\,-{1\over 2}\,
\int_{(0,t]}
X^{(\delta)}_s\, m(ds)\,\Big)\Big)
=
\phi(t)^{\delta/2}
\exp\Big(\,{x\over 2}\,\phi^{\,\prime}(0)\Big)
\end{eqnarray*}
where $v\mapsto\phi(v)$ is a unique positive, continuous and non-increasing function
satisfying
for $0\leq u\leq v\leq t$
\begin{equation}
\label{phi}
\phi^\prime(v)-\phi^{\,\prime}(u)=\int_{(u,v]}\phi(s)\, m(ds),\qquad
\phi(0)=1,
\end{equation}
and which is a constant for $v\geq t$.
By
$\phi^{\,\prime}$ 
we mean the right hand side derivative of $\phi.$
Notice from (\ref{phi}) that $\phi$ is convex in $(0,\infty)$.
Choosing
$$
m(A)=\int_A {2k\,ds\over{(a+s)^2}}+2\gamma\,\varepsilon_{\{1\}}(A),
$$
where $A$ is a Borel set in $(0,\infty)$ and $\varepsilon_{\{1\}}$ is the
Dirac measure at $1,$ yields
\begin{eqnarray}
\label{f12}
\nonumber
&&
\E_0\Big(\exp\Big(\,-{1\over 2}\,
\int_{(0,1]}
X^{(2)}_s\, m(ds)\,\Big)\Big)\\
\nonumber
&&
\hskip2cm
=
\E_0\Big(\exp\Big(\,-k\,
\int_0^1
{X^{(2)}_s\over
{(a+s)^{2}}}\,ds\, -
\gamma\,X^{(2)}_1\Big)\Big)\\
&&
\hskip2cm
=
\phi(1).
\end{eqnarray}
To find the function $\phi,$ we proceed in a similar manner as in
\cite{pitmanyor82} Example 1 p. 432. It follows from (\ref{phi})
that $\phi(v)$ for $v\leq t$ is the (continuous) solution of
\begin{equation}
\label{ode}
\phi''(v)={2k\over{(a+v)^2}}\, \phi(v), 
\end{equation}
satisfying the conditions (notice that $\phi'(1)=0$)
\begin{equation}
\label{odec}
\phi(0)=1,\quad \phi'(1-)=-2\gamma \phi(1).
\end{equation}
It is elementary to check that $y(v)=(a+v)^\alpha$
is a solution of (\ref{ode}) if and only if
$\alpha(\alpha-1)=2k,$ that is
$$
\alpha={1\over 2}(1\pm\sqrt{8k+1})=:\alpha_\pm.
$$
Introducing
$$
y_+(v):=(1+\frac{v}{a})^{\alpha_+}\quad {\rm and}\quad
y_-(v):=(1+\frac{v}{a})^{\alpha_-},
$$
our task is to find constants $A$ and $B$ such that
$$
\phi(v):=A\,y_+(v)+ B\, y_-(v),
$$
fullfills (\ref{odec}).
We skip the detailed computations
and state only the result needed in (\ref{f12}):
$$
\phi(1)= w\,\Big(y_+^{\,\prime}(1)+2\gamma y_+(1)
-(y_-^{\,\prime}(1)+2\gamma y_-(1))\Big)^{-1},
$$
where  $w=\sqrt{8k+1}/a$ is the Wronskian.
The desired formula (\ref{F}) in Proposition \ref{T2} results now
from (\ref{variables2}) and (\ref{f12}) by recalling the definition of
$\gamma$ in (\ref{gamma}) and substituting $r=\log((a+1)/a).$
\end{proof}
\begin{remark}
It is seen from the proof above (take $\gamma=0$) that for $a>0$
$$
\E_0\left(\exp\left(\,-k\,
\int_0^1
{X^{(\delta)}_s\over
{(a+s)^{2}}}\,ds\right)\right)
=\left( 
\frac{\sqrt{8k+1}}{a}\right)^{\delta/2}\ 
\left(y_+^{\,\prime}(1)-y_-^{\,\prime}(1)\right)^{-\delta/2}.
$$
This formula is valid also for $0<\delta<2$ when the boundary point 0
is taken to be reflecting. In particular, $X^{(1)}$ is a Brownian
motion squared. Similarly, for $a>1$ 
\begin{eqnarray}
\label{mansuy}
\nonumber
&&\E_0\left(\exp\left(\,-k\,
\int_0^1
{X^{(\delta)}_s\over
{(a-s)^{2}}}\,ds\right)\right)\\
&&\hskip3cm
=
\left(\frac{\sqrt{8k+1}}{a}\right)^{\delta/2}\ 
\left(x_-^{\,\prime}(1)-x_+^{\,\prime}(1)\right)^{-\delta/2},
\end{eqnarray}
where for $v<a$
$$
x_+(v):=(1-\frac{v}{a})^{\alpha_+}\quad {\rm and}\quad
x_-(v):=(1-\frac{v}{a})^{\alpha_-}.
$$
The formula (\ref{mansuy}) is derived in Mansuy \cite{mansuy03a} 
for squared Brownian motion using different techniques 
(but it is also indicated therein that the result 
can be obtained in the way presented above). See also 
Mansuy \cite{mansuy03b} for  closely related results.
\end{remark}

\subsection{An identity due to Biane and Imhof}
\label{3.5}

Next we consider the identity (\ref{i}) 
(renumbered (\ref{imhof}) below) found by Biane \cite{biane85} and
Imhof \cite{imhof86}. 
This identity is also
observed in \cite{salminennorros01} in connection with a storage
process. The distribution of the random variable $H_\lambda$
featuring on the right hand side of (\ref{imhof}) is in this
context called the RBrownian motion-equilibrium-time-to-emptiness distribution,
see Abate and Whitt \cite{abatewhitt87I}.

We give two proofs
of the Biane--Imhof identity (\ref{imhof}). 
The first one is based on the
Ray--Knight Theorem \ref{r-k2} and the random time change
techniques. Given these tools the proof itself is very short. In
fact, we repeat the idea of the proof of (\ref{343}) in
Proposition \ref{trandufresne}. The second proof is also based on
random time changes but now we work via Tanaka's formula.
This latter presentation is close to the one in \cite{doneygrey89} (see also
Biane \cite{biane85}), but we wish to give it anyway to demonstrate 
the connections between occupation and hitting times.

\begin{proposition} 
{\sl For $\mu>0$
\begin{equation}
\label{imhof} \int_0^\infty  {\bf 1}_{\{
B_s^{(\mu)}\leq 0\}}\, ds \quad{\mathop=^{\rm{(d)}}}\quad
H_\lambda(B^{(\mu)}),
\end{equation}
where $B^{(\mu)}$ is started at 0 and $\lambda$ is an
exponentially (with parameter $2\mu$) distributed random variable
independent of $B^{(\mu)}.$
}
\end{proposition}
{\sl Proof 1 (based on the Ray--Knight Theorem \ref{r-k2}.)}\quad
We use the notation and the structure of the proof of (\ref{343})
in Proposition \ref{trandufresne}. Firstly, by the occupation time
formula and Theorem \ref{r-k2}
$$
\int_0^\infty {\bf 1}_{\{B^{(\mu)}_s<0\}} ds =\int_{-\infty}^0
L^{y}_\infty(B^{(\mu)})\, dy \quad{\mathop=^{\rm{(d)}}}\quad
\int_{0}^\infty Z_{y}\, dy,
$$
where $Z$ satisfies
\begin{equation} \label{371}
    dZ_y= 2\sqrt{Z_y}\,dB_y -2\mu Z_y\,dy.
\end{equation}
The distribution of $Z_0$ is the distribution of
$L^0_\infty(B^{(\mu)})$, i.e., the exponential distribution with
parameter $\mu.$ For the random time change, introduce
$$
C_y:=\int_0^y Z_u\, du
$$
and let $c$ denote the inverse of $C.$ It follows from
(\ref{371})
$$
Z_{c_y}-Z_0 =2\,\beta_{y}-2\mu\,y,
$$
where $\beta$ is a Brownian motion. Consequently, because $Z_y\to 0$ as
$y\to\infty$ we obtain (\ref{imhof}).~$\square$
{\sl Proof 2 (based on the Tanaka formula).}\quad Consider 
$$
\bigl(B_t^{(\mu)}\bigr)^-=
-\int_0^t  {\bf 1}_{\{
B_s^{(\mu)}\leq 0\}}\, dB_s^{(\mu)}
+{1\over 2}\ L^0_t(B^{(\mu)}),
$$
where $L^0(B^{(\mu)})$ is the local time of
$B^{(\mu)}$ at 0 (with respect to the Lebesgue measure).
Introduce
$$
A^-_t:=\int_0^t  {\bf 1}_{\{
B_s^{(\mu)}\leq 0\}}\, ds,\quad {\rm and}\quad \alpha^-_t:=
\inf\{s:\ A_s^->t\}.
$$
Then letting $\ell^{(\mu)}_t= {1\over 2}\ L^0_t(B^{(\mu)})$
\begin{eqnarray*}
\bigl(B_{\alpha^-_t}^{(\mu)}\bigr)^-=-B_{\alpha^-_t}^{(\mu)}&=&
-\int_0^{\alpha^-_t}  {\bf 1}_{\{
B_s^{(\mu)}\leq 0\}}\, dB_s^{(\mu)}
+\ell^{(\mu)}_{\alpha^-_t}\\
&=&
-\int_0^{\alpha^-_t}  {\bf 1}_{\{
B_s^{(\mu)}\leq 0\}}\, dB_s-
\mu\int_0^{\alpha^-_t}  {\bf 1}_{\{
B_s^{(\mu)}\leq 0\}}\, ds
+\ell^{(\mu)}_{\alpha^-_t}\\
&=&
-\int_0^{\alpha^-_t}  {\bf 1}_{\{
B_s^{(\mu)}\leq 0\}}\, dB_s-
\mu\, t
+\ell^{(\mu)}_{\alpha^-_t}.
\end{eqnarray*}
It is a straightforward application of the result due to Dambis, and
Dubins and Schwarz (see Revuz and Yor \cite{revuzyor01} p. 181) that the
process given by
$$
\beta_t:=\int_0^{\alpha^-_t}  {\bf 1}_{\{
B_s^{(\mu)}\leq 0\}}\, dB_s
$$
is an $\cF_{\alpha^-_t}$--Brownian motion. Notice that
because $B_{\alpha^-_t}^{(\mu)}\leq 0$  we have for all $t\geq 0$
\begin{equation}
\label{smaller}
\ell^{(\mu)}_{\alpha^-_t}\geq \beta_t+\mu\, t.
\end{equation}
Clearly,
$
A^-_t= A^-_{\Lambda_0}
$
for $t\geq \Lambda_0:=\Lambda_0(B^{(\mu)}):=\sup\{ t:\ B^{(\mu)}_t=0\}.$
Consequently, $\{B_{\alpha^-_t}^{(\mu)}:\ t\geq 0\}$ is defined only
for $t<A^-_{\Lambda_0}$ and
$$
\lim_{t\to A^-_{\Lambda_0}} B_{\alpha^-_t}^{(\mu)}= B_{\Lambda_0}^{(\mu)}=
0.
$$
Because $\alpha^-$ is the inverse of $A^-$ and $\ell^{(\mu)}$ does not
increase after $\Lambda_0$ we obtain
$$
0=-\beta_{A^-_{\Lambda_0}}-\mu A^-_{\Lambda_0}+
\ell^{(\mu)}_{\Lambda_0}=
-\beta_{A^-_{\Lambda_0}}-\mu A^-_{\Lambda_0}+
\ell^{(\mu)}_{\infty}.
$$
From (\ref{smaller}) it now follows
that
$$
A^-_\infty=A^-_{\Lambda_0}=\inf\{t:\ \beta_t+\mu t=
\ell^{(\mu)}_{\infty}\}.
$$
Because $\beta$ is a Brownian motion and the
$\P_0$--distribution of $\ell^{(\mu)}_{\infty}$ is
exponential with parameter $2 \mu$ it remains to prove that
$\beta$ and $\ell^{(\mu)}_{\infty}$ are independent. To do this,
let
$$
A^+_t:=\int_0^t  {\bf 1}_{\{
B_s^{(\mu)}\geq 0\}}\, ds,\quad {\rm and}\quad \alpha^+_t:=
\inf\{s:\ A_s^+>t\},
$$
and proceed as above to obtain
\begin{eqnarray*}
\bigl(B_{\alpha^+_t}^{(\mu)}\bigr)^+=B_{\alpha^+_t}^{(\mu)}&=&
\int_0^{\alpha^+_t}  {\bf 1}_{\{
B_s^{(\mu)}\geq 0\}}\, dB_s^{(\mu)}
+\ell^{(\mu)}_{\alpha^+_t}\\
&=&
\gamma_t+\mu t+\ell^{(\mu)}_{\alpha^+_t},
\end{eqnarray*}
where $\{\gamma_t:\ t\geq 0\}$ is a Brownian motion.
Because
$
B_{\alpha^+_t}^{(\mu)}
$
is well defined and non-negative for all
$t\geq 0$ we deduce from Skorokhod's reflection equation that
$$
\ell^{(\mu)}_{\alpha^+_t}=\sup_{0\leq s\leq t}\{-\gamma_s-\mu s\},
$$
and letting $t\to\infty$ gives
$$
\ell^{(\mu)}_\infty=\sup_{s\geq 0}\{-\gamma_s-\mu s\}.
$$
From the extended form of Knight's theorem (see \cite{revuzyor01}
p. 183) we know that $\beta$ and $\gamma$ are
independent and, therefore, also $\ell^{(\mu)}_\infty $ and $\beta$ are
independent.~$\square$

\subsection{LeGall's identity}

From a representation formula in LeGall \cite{legall85} we can
deduce the result in (\ref{legall}) below, see 
Donati-Martin and Yor \cite{donatimartinyor97} p. 1044 and 1052.
In this section we prove (\ref{legall}) using time reversal argument
and time changes. This  nice application of the classical time
reversal result was pointed out to us by Y. Hariya (in the context of 
(\ref{hariya})). 

\begin{proposition}
{\sl  The following formula holds: 
\begin{equation}
\label{legall}
\int_0^\infty \exp\big(-2\,R^{\,(3)}_s\big)\, ds
\ {\mathop=^{\rm{(d)}}}\
H_{1}(R^{\,(2)}),
\end{equation}
where the Bessel processes are started from 0.
}
\end{proposition}
\begin{proof}
For $x>0$ let $\Lambda_x(R^{\,(3)}):=\sup\{t:\ R^{\,(3)}_t =x\}.$
Then 
$$
\int_0^{\Lambda_x(R^{\,(3)})} \exp\big(-2\,R^{\,(3)}_s\big)\, ds 
\ {\mathop=^{\rm{(d)}}}\
\int_0^{H_0(B)} \exp\big(-2\,B_s\big)\, ds,
$$
where $B$ is a standard Brownian motion started from $x.$
Let for $t>0$
$$
A_t:=\int_0^t  \exp\big(-2\,B_s\big)\, ds
$$
and, as usual, $\alpha_t$ is its inverse. By the Lamperti
transformation (\ref{lam1}), the process $\{Z_t:\ t\geq 0\}$ where
$Z_t:=\exp(-B_{\alpha_t})$ is a 2-dimensional Bessel process. Clearly,
$$
Z_0=\exp(-x)\quad {\rm and}\quad 0<Z_t<1\quad \forall\ 
t<A_{H_0}.
$$
Consequently,
$$
A_{H_0}
\ {\mathop=^{\rm{(d)}}}\
H_{1}(Z),
$$
and letting $x\to\infty$ proves (\ref{legall}).
\end{proof}

\subsection{Hariya's identity}

We learned the identity in the next proposition from Y. Hariya 
\cite{hariya04} and offer here a proof which differs from Hariya's
proof and is in the spirit of the present paper. 

\begin{proposition}
{\sl The following identity holds:
\begin{equation}
\label{hariya}
\int_0^\infty (1+R^{\,(3)}_s)^{-2\gamma}\, ds
\ {\mathop=^{\rm{(d)}}}\
H_{1/(\gamma-1)}(R^{\,(\delta)}),
\end{equation}
where $\gamma>1,$ $\delta=(2\gamma-1)/(\gamma-1)$
and the Bessel processes are started from~0.
}
\end{proposition}
\begin{proof}
By the scaling property of Bessel
processes and 
the occupation time formula the right hand side 
of (\ref{hariya}) can be written as 
\begin{eqnarray*}
&&H_{1/(\gamma-1)}(R^{(\delta)})
\ {\mathop=^{\rm{(d)}}}\
(\gamma-1)^{-2}H_1(R^{(\delta)})\\
&&\hskip2.5cm\ {\mathop=^{\rm{(d)}}}\
\frac{1}{(\gamma-1)^2}
\int_0^{1}
L^y_{H_{1}}(R^{(\delta)})\, dy.
\end{eqnarray*}
Consequently, by Theorem \ref{r-k3} in Appendix  
(with the notations as therein),
(\ref{hariya}) is now equivalent with 
\begin{equation}
\label{hariya2}
\int_0^\infty \frac{Z^{(2)}_y}{(1+y)^{2\gamma}}\, dy
\ {\mathop=^{\rm{(d)}}}\
(\delta-2)\,\int_0^1\frac{\widehat
Z^{(2)}_{y^{\delta-2}}}{y^{\delta-3}}\,dy.
\end{equation}
In particular, when $\gamma=2$ (and $\delta =3$) the identity
(\ref{hariya2}) takes the simple form
\begin{equation}
\label{hariya3}
\int_0^\infty \frac{Z^{(2)}_y}{(1+y)^{4}}\, dy
\ {\mathop=^{\rm{(d)}}}\
\int_0^1 \widehat Z^{(2)}_{y}\,dy.
\end{equation}
To prove (\ref{hariya}) via (\ref{hariya2}) we use the well known fact that
\begin{equation}
\label{bridge}
\{ Z^{(2)}_s:\ s\geq 0\}
\ {\mathop=^{\rm{(d)}}}\
\{ (1+s)^2\widehat Z^{(2)}_{1/(1+s)}:\ s\geq 0\}.
\end{equation}
Notice that (\ref{hariya3}) follows directly from (\ref{bridge}).
However, to obtain (\ref{hariya2}) we have to work little more.
Substituting $x=y^{\delta-2}$ on the right hand side of (\ref{hariya2})
shows that (\ref{hariya2}) is equivalent with
\begin{equation}
\label{hariya4}
\int_0^\infty \frac{Z^{(2)}_y}{(1+y)^{2\gamma}}\, dy
\ {\mathop=^{\rm{(d)}}}\
\int_0^1 x^{2(\gamma-2)}\widehat Z^{(2)}_{x}\,dx.
\end{equation}
By the representation (\ref{bridge}) the left hand side of
(\ref{hariya4}) takes the form
$$
\int_0^\infty \frac{Z^{(2)}_y}{(1+y)^{2\gamma}}\, dy
\ {\mathop=^{\rm{(d)}}}\
\int_0^\infty \frac{\widehat Z^{(2)}_{1/(1+y)}}{(1+y)^{2\gamma-2}}\, dy,
$$
and substituting here $x=1/(1+y)$ leads to the right hand side of
(\ref{hariya4}), completing the proof.
\end{proof}

\begin{remark}
\label{Ciesielski--Taylor}
Consider the Ciesielski--Taylor identity (\ref{c-t}):
\begin{equation}
\label{c-t2}
\int_0^\infty {\bf 1}_{\{R^{(\delta+2)}_s< 1\}}\, ds
\quad{\mathop=^{\rm{(d)}}}\quad
H_1(R^{(\delta)}),
\end{equation}
where $R^{(\delta)}$ is a Bessel process of dimension $\delta>0$
started at 0. It is possible to express (\ref{c-t2}) in alternative
forms using Hariya's identity and some simple transformations.
Indeed, by Theorem \ref{r-k2} (b) we rewrite (\ref{c-t2}) first in
the form (see Yor \cite{yor92c} p. 52)
\begin{equation}
\label{c-t3}
\frac{1}{\delta-2}\int_0^1 \frac{\widehat Z^{(2)}_{y^{\delta-2}}}{\
y^{\delta-3}\ }\, dy
\ {\mathop=^{\rm{(d)}}}\
\frac{1}{\delta}\int_0^1 \frac{Z^{(2)}_{y^{\delta}}}{\ y^{\delta-1}\ }\, dy.
\end{equation}
Applying (\ref{hariya2}) on the left hand side of (\ref{c-t3}) yields
\begin{equation}
\label{c-t4}
\frac{1}{(\delta-2)^2}\int_0^\infty \frac{ Z^{(2)}_{y}}{(1+y)^{2\gamma}}\, dy
\ {\mathop=^{\rm{(d)}}}\
\frac{1}{\delta}\int_0^1 \frac{Z^{(2)}_{y^{\delta}}}{\ y^{\delta-1}\ }\, dy.
\end{equation}
Making the change of variables $x=y^\delta$ on the right hand side of
(\ref{c-t4}) and recalling that
$\gamma=(\delta-1)/(\delta-2)$ leads to
\begin{equation}
\label{c-t5}
\frac{1}{(\delta-2)^2}\int_0^\infty \frac{ Z^{(2)}_{y}}{(1+y)^{2(\delta-1)/(\delta-2)}}\, dy
\ {\mathop=^{\rm{(d)}}}\
\frac{1}{\delta^2}\int_0^1 \frac{Z^{(2)}_{x}}{x^{2(\delta-1)/\delta}}\, dx.
\end{equation}
Further, substituting on the right hand side of (\ref{c-t5}) $y=1/x$
and using the time inversion property of Bessel processes we obtain
\begin{eqnarray*}
&&\frac{1}{(\delta-2)^2}\int_0^\infty \frac{ Z^{(2)}_{y}}
{(1+y)^{2(\delta-1)/(\delta-2)}}\, dy\\
&&\hskip4cm\ {\mathop=^{\rm{(d)}}}\
\frac{1}{\delta^2}\int_1^\infty y^{2/\delta}\ Z^{(2)}_{1/y}\, dy\\
&&\hskip4cm\ {\mathop=^{\rm{(d)}}}\
\frac{1}{\delta^2}\int_1^\infty \frac{ Z^{(2)}_{z}}
{z^{2(\delta-1)/\delta}}\, dz\\
&&\hskip4cm\ {\mathop=^{\rm{(d)}}}\
\frac{1}{\delta^2}\int_0^\infty \frac{ Y^{(2)}_{u}}
{(1+u)^{2(\delta-1)/\delta}}\, du,
\end{eqnarray*}
where $Y^{(2)}$ denotes the process $Z^{(2)}$ started with an
exponential distribution with parameter $1/2$ (which is the
distribution of $Z^{(2)}_1$ when started from 0).
\end{remark}

\section{Feynman-Kac approach to perpetual integral functionals}
\label{sec5}

In this section we show how the Feynman-Kac formula can be used to
find the Laplace transform of a perpetual integral functional
$I_\infty(f)$ where $f$ satisfies the integrability condition (\ref{K1}):
$$
\int^\infty f(y)dy<\infty.
$$
We remark also
that in special cases one can find the law of a perpetual functional
by limiting procedures but for a general
characterization the problem has to be analyzed more carefully.

For treatments of Feynman-Kac formula, we refer to Durrett
\cite{durrett84} and Karatzas and Shreve \cite{karatzasshreve91}. See
also Jeanblanc, Pitman and Yor \cite{jeanblancpitmanyor97} for
connections with excursion theory.

Consider for $\gamma>0$ (and $\mu>0$)
the function
$$
x\mapsto\Psi_\gamma(x)=\E_x\Big(\exp\bigr(-I_\infty(\gamma f)\big)\Big)=
\E_x\Bigr(\exp\Bigr(-\gamma\int_0^\infty f(B^{(\mu)}_s)\, ds
\Bigr)\Bigr).
$$
By the simple Markov property it is seen that the process
$$
\Big\{\Psi_\gamma(B^{(\mu)}_t)\,
\exp\Big(-\gamma\int_0^t f(B^{(\mu)}_s)\, ds
\Big):\, t\geq 0\Big\}
$$
is a bounded martingale.

\begin{proposition}
\label{Psi}
{\sl
The function $\Psi_\gamma$ is non-decreasing and satisfies
 for all $x$ and $t\geq 0$
\begin{equation}
\label{invariant}
\Psi_\gamma(x)=\E_x\Big(
\Psi_\gamma(B^{(\mu)}_t)\ \exp\Big(-\gamma\int_0^t f(B^{(\mu)}_s)\, ds
\Big)\Bigr).
\end{equation}
Moreover,
\begin{equation}
\label{limit1}
\lim_{x\to\infty}\Psi_\gamma(x)=1.
\end{equation}
}
\end{proposition}
\begin{proof}
The formula (\ref{invariant}) is immediate from the martingale
property. Let $x<y$ and apply the optional stopping theorem to obtain
\begin{eqnarray}
\label{b0}
\Psi_\gamma(x)&=&
\E_x\Big(
\Psi_\gamma(B^{(\mu)}_{H_y})\,
\exp\Big(-\gamma\int_0^{H_y} f(B^{(\mu)}_s)\, ds
\Big)\Big)\\
\nonumber
&\leq &
\Psi_\gamma(y),
\end{eqnarray}
which shows that $\Psi_\gamma$ is non-decreasing.
Next notice that
$$
\Psi_\gamma(B^{(\mu)}_t)\,
\exp\Big(-\gamma\int_0^t f(B^{(\mu)}_s)\, ds
\Big)
=
\E\Big(\exp\Big(-\gamma\int_0^\infty f(B^{(\mu)}_s)\, ds
\Big)\,|\,\cF_t\Big).
$$
Consequently, by the martingale convergence theorem,
$$
\Psi_\gamma(\infty)\, \exp\Big(-\gamma\int_0^\infty
f(B^{(\mu)}_s)\, ds \Big) = \exp\Big(-\gamma\int_0^\infty
f(B^{(\mu)}_s)\, ds\Big)
$$
showing (\ref{limit1}).
\end{proof}
\begin{remark}
\label{remark1}
Notice from (\ref{b0}) that for $x\leq y$
$$
\E_x\Big(
\exp\Bigr(-\gamma\int_0^{H_y} f(B^{(\mu)}_s)\, ds
\Bigr)\Big)
=\frac{\Psi_\gamma(x)}{\Psi_\gamma(y)}.
$$
\end{remark}

In many cases the function $\Psi_\gamma$ can be found by solving a
second order ODE. This is formulated in the following
\begin{proposition}
\label{ODE}
{\sl
Assume that $f$ is piecewise continuous and satisfies the
integrability condition (\ref{K1}).
Then $x\mapsto\Psi_\gamma(x)$ is the unique positive,
non-decreasing and continuously differentiable solution of
the problem
\begin{eqnarray}
\label{system}
&&
{1\over 2}\, v''(x)+\mu\,v'(x)-\gamma\, f(x)\,v(x)=0,\\
\nonumber
&&
\lim_{x\to\infty}v(x)=1.
\end{eqnarray}
}
\end{proposition}
\begin{proof}
Notice that
$B^{(\mu)}$ killed according to the additive functional
$$
I_t(\gamma f)=
\gamma\int_0^t f(B^{(\mu)}_s)\, ds
$$
is a linear diffusion $B^\bullet$ (in the sense of Ito and McKean
\cite{itomckean74}), and its basic characteristics (speed measure,
scale function and killing measure) can be determined explicitly
(see \cite{borodinsalminen02} No. II.9 p. 17). From
(\ref{invariant}) it follows that $\Psi_\gamma$ is an invariant
function for $B^\bullet.$ It is well known that invariant
functions of a linear diffusion are continuous (see Dynkin
\cite{dynkin65} Vol. II p. 7), and differentiable when the scale
function is differentiable (see Salminen \cite{salminen85} p. 93).
 From the representation
theory of excessive functions we know that for $B^\bullet$ there
exist two invariant functions, denoted $\varphi_\gamma$ and
$\psi_\gamma$ and  called {\sl fundamental solutions} or {\sl
extreme invariant functions}, such that if $h$ is an arbitrary
invariant function then there exists constants $c_1\geq 0$ and
$c_2\geq 0$ such that $h=c_1\, \psi_\gamma + c_2\,
\varphi_\gamma.$ Moreover, we have 
$$
\P_x(H_z(B^\bullet)<\infty)=\cases{
{\displaystyle\frac{\varphi_\gamma(x)}{\varphi_\gamma(z)}},& $x\geq
  z,$\cr
&\cr {\displaystyle\frac{\psi_\gamma(x)}{\psi_\gamma(z)}},& $x\leq
z,$\cr}.
$$
From this representation it follows that $\varphi$ and $\psi$ 
solve (\ref{system}) on the
intervals of continuity of $f.$ 
Next notice that, because $\mu>0,$
$$
\lim_{z\to -\infty}\P_x(H_z(B^\bullet)<\infty)=0
$$
giving $\varphi_\gamma(z)\to +\infty$ as $z\to -\infty.$ Consequently,
all invariant non-decreasing functions are multiples of $\psi_\gamma.$
In particular,  $\Psi_\gamma$ is a multiple of  $\psi_\gamma$ and,
hence, the condition (\ref{limit1}) determines  $\Psi_\gamma$
uniquely.
\end{proof}
\begin{remark} At the points of discontinuity of $f$ the function
$\Psi_\gamma$ usually fails to be two times differentiable.
\end{remark}

\begin{example}
\label{yor1}
We compute the Laplace
transform of the functional
$$
\int_0^\infty
{{(a+\exp(B^{(1/2)}_s))^{-2}}}\, ds
$$
appearing in (\ref{y}). Consider, for a moment, the case with general  
$\mu>0.$ Taking 
$$
f(x)=(a+\exp({x}))^{-2}
$$
in (\ref{system}) gives us the equation 
$$
{1\over 2}\, v''(x)+\mu\,v'(x)-\gamma\, (a+\exp({x}))^{-2}\,v(x)=0.
$$
Putting here $x=\ln y$ and $g(y)=v(\ln y)$ yields 
\begin{equation}
\label{de01}
{1\over 2}\,y^2\, g''(y)+(\mu+{1\over 2})\,y\,g'(y)
-\gamma\ (a+y)^{-2}\,g(y)=0,
\end{equation}
which is, of course, the corresponding equation for geometric
Brownian motion. By
Kamke \cite{kamke43}
2.394 p. 497 this equation can be solved for $\mu =1/2$
by making the substitution
$$
\eta(\xi)=g(y),\quad \xi={\sqrt{2\gamma}\over a}\ln
\bigl({y\over{y+a}}\bigr),
$$
which transforms (\ref{de01}) to the following 
\begin{equation}
\label{de1}
\sqrt{2\gamma}\ \eta''+a\ \eta'=\sqrt{2\gamma}\ \eta.
\end{equation}
Letting $\beta:=a/2\sqrt{2\gamma}$ the general solution of (\ref{de1})
can be written as
$$
\eta(\xi)=A\,\exp\Big(-(\sqrt{1+\beta^2}+\beta)\xi\Big)
+B\,\exp\Big((\sqrt{1+\beta^2}-\beta)\xi\Big).
$$
Consequently, the increasing solution of (\ref{de01}) is
\begin{eqnarray*}
\psi(y)&=&\exp\Bigl((\sqrt{1+\beta^2}-\beta){1\over{2\beta}}
\ln\bigl({y\over{a+y}}\bigr)\Bigr)\\
&=&\Bigl({y\over{a+y}}\Bigr)^{2^{-1}(\sqrt{1+\beta^{-2}}-1)}.
\end{eqnarray*}
Notice that $\psi(\infty)=1,$ and it follows 
\begin{eqnarray}
\label{01}
\nonumber
&&\E_x\Bigr(\exp\Bigr(-\gamma\int_0^\infty
{{(a+\exp(B^{(1/2)}_s))^{-2}}}\, ds
\Bigr)\Bigr)\\
&&\hskip3cm
=
\Bigl({\exp(x)\over{a+\exp(x)}}
\Bigr)^{(2\,a)^{-1}\,(\sqrt{a^2+8\gamma}-a)}.
\end{eqnarray}
For a geometric Brownian motion $X$ with $X_0=x>0$ defined via 
$$
X_s=\exp(B^{(1/2)}_s),\quad  B^{(1/2)}_0=x,
$$ 
the formula (\ref{01}) takes the form
\begin{equation}
\label{1}
\E_x\Bigl(\exp\Bigl(-\gamma\int_0^{\infty} {1\over{(a+X_s)^{2}}}\, ds
\Bigr)\Big)
=\Bigl({x\over{a+x}}\Bigr)^{(2\,a)^{-1}\,(\sqrt{a^2+8\gamma}-a)}.
\end{equation}
The identity in law in (\ref{y}) can be deduced from (\ref{01}) (or (\ref{1})).
Notice also that substituting in (\ref{01}) $x=0$ and 
letting $a\to 0$ we obtain
a special case of the identity (\ref{d-y}):
$$
\E_0\Bigl(\exp\Bigl(-\gamma\int_0^{\infty}   e^{2B_s-s}\, ds
\Bigr)\Big)
=e^{-\sqrt{2\gamma}},
$$
i.e.,
$$
\int_0^\infty
 e^{2B_s-s}\, ds\ {\mathop=^{\rm{(d)}}}\
\Bigl(2\,Z_{1/2}\Bigr)^{-1},
$$
where $Z_{1/2}$ is a $\Gamma$--distributed r.v. with parameter $1/2.$
\end{example}

Using Propositions \ref{Psi} and \ref{ODE} we derive an
interesting result due to Biane \cite{biane85} which characterizes the
law of a perpetual integral functional  of $B^{(\mu)},\
\mu>0,$ restricted on $\R_-$ in terms of the same but unrestricted
functional of another diffusion stopped at the first hitting time.
We remark that in \cite{biane85} a more general situation (not only
$B^{(\mu)}$) is considered. However, the main interest in
\cite{biane85} is focused on occupation times, the aim being to
generalize the Ciesielski--Taylor identity (\ref{c-t}). The result
in our Proposition \ref{biane} is extracted from
Remarque p. 295 in \cite{biane85} and formulated
for $B^{(\mu)}.$
\begin{proposition}
\label{biane} {\sl Let $f$ be a positive $\cC^1$-function such
that
\begin {equation}
\label{S}
\int_{-\infty}f(y)\, {\rm e}^{-2\mu y}\ dy =\infty,
\end{equation}
and $X$ a diffusion with the generator
$$
{1\over 2}{d^2\over{dx^2}}
+\Big(\mu -{1\over {2}}\,{f'(x)\over{f(x)}}\Big)\,{d\over{dx}}.
$$
Then $H_0(X)<\infty$ a.s. if $X_0<0,$  and, moreover,
$$
I^-_\infty(f):=\int_0^\infty f(B^{(\mu)}_s)\, {\bf
1}_{\{B^{(\mu)}_s<0\}}\, ds \quad{\mathop=^{\rm{(d)}}}\quad
\int_0^{H_0(X)} f(X_s)\, ds
$$
where $X_0$ is taken to be exponentially  distributed on $(-\infty,0)$
with parameter $2\mu.$
}
\end{proposition}
\begin{proof}
Notice that the condition (\ref{S}) means that the scale function
$S^X$ of $X$ given by (cf. \cite{borodinsalminen02} II.9
p. 17)
$$
S^X(x)=\int^x f(y)\,{\rm e}^{-2\mu y}\, dy
$$
satisfies
$$
\lim_{x\to -\infty} S^X(x)=-\infty.
$$
This implies $H_0(X)< +\infty$ a.s. when $X_0<0.$ From
Proposition \ref{ODE} we know 
$$
\Psi_\gamma(x):=\E_x\Big(\exp\Big(-\gamma I^-_\infty(f)\Big)\Big)
$$
is the unique, non-decreasing function such that 
\begin{equation}
\label{b1} {1\over 2}\, \Psi''(x)+\mu\,
\Psi'(x)=\gamma{\bf 1}_{(-\infty,0)}(x)\,f(x)\, \Psi(x)
\end{equation}
and $\lim_{x\to +\infty} \Psi(x)=1.$ For $x>0$ we clearly have
$$
\Psi_\gamma(x)=\P_x(H_0(B^{(\mu)})=+\infty)
+\P_x(H_0(B^{(\mu)})<+\infty)\,\Psi_\gamma(0),
$$
and, hence, it is enough to compute $\Psi_\gamma(x)$ for $x\leq
0.$ For this, consider the equation
\begin{equation}
\label{b2}
{1\over 2}\, u''(x)+\mu\, u'(x)=\gamma \,f(x)\, u(x).
\end{equation}
Let $\psi_\gamma$ and $\varphi_\gamma$ denote the fundamental
non-decreasing and non-increasing, respectively, solutions of
(\ref{b1}), and, similarly, $\widehat \psi_\gamma$ and $\widehat
\varphi_\gamma$ are the fundamental non-decreasing and
non-increasing, respectively, solutions of (\ref{b2}). Notice that
$f$ does not have to satisfy the integrability condition
(\ref{K1}). However, for  (\ref{b2}), we can still argue as in the
proof of Proposition \ref{ODE} that all invariant non-decreasing
functions are multiples of $\widehat \psi_\gamma.$ Using
continuity and differentiability requirements, $\psi_\gamma$ can
be expressed in terms of $\widehat \psi_\gamma$ as follows
$$
\psi_\gamma(x)=\cases{\widehat\psi_\gamma(x),&$x\leq 0,$\cr
&\cr
 S(x)\,{\displaystyle\frac{\widehat\psi^{\,\prime}_\gamma(0)}{S^{\,\prime}(0)}}
+\widehat\psi_\gamma(0)
-{\displaystyle\frac{S(0)\,\widehat\psi^{\,\prime}_\gamma(0)}
{S^{\,\prime}(0)},}&$x\geq 0,$\cr}
$$
where $S(x)=-\exp(-2\mu\,x)$ is the scale function of  $B^{(\mu)}.$
Consequently, for $x\leq 0$
$$
\Psi_\gamma(x)={\displaystyle\frac{\psi_\gamma(x)}
{\psi_\gamma(+\infty)}}=
{\displaystyle\frac{2\mu\,\widehat\psi_\gamma(x)}
{2\mu\,\widehat \psi_\gamma(0)+\widehat\psi^{\,\prime}_\gamma(0)}},
$$
and, in particular, for $x=0$
\begin{equation}
\label{b21}
\E_0\Big(\exp\Big(-\gamma\,\int_0^\infty f(B^{(\mu)}_s)\, {\bf
1}_{\{B^{(\mu)}_s<0\}}\, ds\Big)\Big)=
{\displaystyle\frac{2\mu\,\widehat\psi_\gamma(0)}
{2\mu\,\widehat \psi_\gamma(0)+\widehat\psi^{\,\prime}_\gamma(0)}}.
\end{equation}
To proceed, define for $y\leq 0$
$$
\widetilde \psi_\gamma(y):=
2\mu\,\widehat \psi_\gamma(y)+\widehat\psi^{\,\prime}_\gamma(y),
$$
and notice that $\widetilde \psi^{\,\prime\prime}_\gamma$ exists and
is continuous because $f\in\cC^1.$
Using the fact that $\widehat \psi_\gamma$ solves (\ref{b2})
it is straightforward to verify that $\widetilde \psi^{\,
\prime}_\gamma>0,$ i.e., $\widetilde \psi_\gamma$ is increasing,
and that $\widetilde \psi_\gamma$ is a solution of the ODE
$$
{1\over 2}u''(x) +\Big(\mu -{1\over
{2}}\,{f'(x)\over{f(x)}}\Big)\,u'(x) =\gamma\,f(x)\, u(x),\quad
x\leq 0.
$$
By It\^o's formula, the process
$$
\Big\{\widetilde \psi_\gamma(X_{t\wedge H_0(X)})\,
\exp\bigl(-\gamma\int_0^{t\wedge H_0(X)} f(X_s)\, ds \bigr):\, t\geq
0\Big\}
$$
is a martingale and, further, because it is bounded, we obtain for
$X_0=x<0$ by the dominated convergence theorem
\begin{equation}
\label{b3}
 \E_x\Big(\exp\Big(-\gamma\int_0^{H_0(X)} f(X_s)\, ds
\Big)\Big) ={\displaystyle \frac{\widetilde \psi_\gamma(x)}
{\widetilde \psi_\gamma(0)}}.
\end{equation}
Observe that
$$
\int_{-\infty}^0 \widetilde \psi_\gamma(x)\, 2\mu\, {\rm e}^{2\mu
x}\, dx = 2\mu \int_{-\infty}^0 \big(2\mu \widehat \psi_\gamma(x)
+ \widehat \psi^{\,\prime}_\gamma(x)\big)\, {\rm e}^{2\mu x}\, dx
= 2\mu\widehat \psi_\gamma(0).
$$
Consequently, if $X_0$ is exponentially distributed on
$(-\infty,0)$ with parameter $2\mu,$ (\ref{b21}) and
(\ref{b3}) lead to
\begin{eqnarray*}
&&
\int_{-\infty}^0\, 2\mu\, {\rm e}^{2\mu x}\,
\E_x\Big(\exp\Big(-\gamma\int_0^{H_0(X)} f(X_s)\, ds \Big)\Big)\, dx\\
&&\hskip3cm =\E_0\Big(\exp\Big(
-\gamma \int_0^\infty f(B^{(\mu)}_s)\, {\bf
1}_{\{B^{(\mu)}_s<0\}}\, ds\Big)\Big),
\end{eqnarray*}
as claimed.
\end{proof}

\begin{example} Consider the functional 
$$
I:=\int_0^\infty \exp(-2\,B^{(\mu)}_s)\, {\bf
1}_{\{B^{(\mu)}_s<0\}}\, ds.
$$
Recall from (\ref{333}), that there exists a Bessel process 
$R^{(2-2\mu)}$ started from $1$ such that
$$
I\,=\,
\int_0^\infty
{\bf 1}_{\{R^{(2-2\mu)}_{s}> 1\}}\, ds\quad{\rm a.s.}
$$
As an application of Proposition \ref{biane} we
derive a new charaterization of the distribution of $I.$
Taking $\displaystyle f(x)={\rm e}^{-2x}$ it is seen that the diffusion
$X$ in Proposition \ref{biane} is a Brownian motion 
with drift $\mu+1.$ Consequently,
$$
I
\quad{\mathop=^{\rm{(d)}}}\quad
\int_0^{H_0(B^{(\mu+1)})} \exp\left(-2\,B^{(\mu+1)}_s\right)\, ds,
$$
where $B^{(\mu+1)}_0$ is exponentially distributed on $(-\infty,0)$ 
with parameter $2\mu.$ To develop this further, let $x>0$ 
and assume that $B^{(\mu+1)}_0=-x<0.$ We have 
\begin{eqnarray*}
&&\int_0^{H_0(B^{(\mu+1)})} \exp\left(-2\,B^{(\mu+1)}_s\right)\, ds
\\
&&
\hskip3cm
 \quad{\mathop=^{\rm{(d)}}}\quad
\int_0^{H_x(\hat B^{(\mu+1)})} \exp\left(-2\,(\hat B^{(\mu+1)}_s-x)
\right)\, ds
\\
&&\hskip3.45cm 
=\ 
{\rm e}^{2x}\
\int_0^{H_x(\hat B^{(\mu+1)})} \exp\left(-2\,\hat B^{(\mu+1)}_s
\right)\, ds
\\
&&\hskip3cm 
\quad{\mathop=^{\rm{(d)}}}\quad
{\rm e}^{2x}\
\inf\{t\,:\, R^{(-2\mu)}_t={\rm e}^{-x}\},
\end{eqnarray*}
where $\hat B^{(\mu+1)}$ is a Brownian motion with drift $\mu+1$
started from 0 and in the last step Remark \ref{stoppedfunctional}
is applied. Using the scaling property of Bessel processes we obtain
\begin{eqnarray*}
&&\inf\{t\,:\, R^{(-2\mu)}_t={\rm e}^{-x}\}
= 
\inf\{t\,:\,{\rm e}^{x}\, R^{(-2\mu)}_t=1\}
\\
&&\hskip4cm
=
\inf\{{\rm e}^{-2x}\,t\,:\,{\rm e}^{x}\, 
R^{(-2\mu)}_{{\rm e}^{-2x}\,t}=1\}
\\
&&
\hskip3.6cm
\quad{\mathop=^{\rm{(d)}}}\quad
{\rm e}^{-2x}\,\inf\{t\,:\, \hat R^{(-2\mu)}_t=1\},
\end{eqnarray*}
where the Bessel process $\hat R^{(-2\mu)}$ is started from 
${\rm  e}^{x}.$ Consequently, 
\begin{equation}
\label{b11}
\int_0^\infty \exp(-2\,B^{(\mu)}_s)\, {\bf
1}_{\{B^{(\mu)}_s<0\}}\, ds.
\quad{\mathop=^{\rm{(d)}}}\quad
H_1(\hat R^{(-2\mu)}),
\end{equation}
where $\hat R^{(-2\mu)}_0$ is distributed as ${\rm e}^\xi$ with $\xi$ 
exponentially distributed with parameter $2\mu.$ Elementary
computations show that 
\begin{equation}
\label{speed}
\P(\hat R^{(-2\mu)}_0>z)= z^{-2\mu},\quad z\geq 1.
\end{equation}
It is interesting to notice that the right hand side of (\ref{speed}) 
when extended to a measure on the whole of 
$\R_+$ can be viewed as the speed measure of 
$\hat R^{(-2\mu)}$ (see, e.g., \cite{borodinsalminen02} A1.21 p. 133).
\end{example}

\section{Appendix on Ray--Knight theorems}
\label{sec6}
For an easy reference, we recall here the Ray--Knight theorems used
in this paper (see  Yor \cite{yor92c} and
\cite{borodinsalminen02} for more complete statements).

\begin{theorem}
\label{r-k1}
{\sl Let $B$ be a standard Brownian motion started from 0 and
  $L^y_{H_1(B)}(B)$ its local time (with respect to the Lebesgue
measure) at level $y\leq 1$ up to $H_1(B).$
Then the local time process $\{L^{1-y}_{H_1}(B):\ y\geq 0\}$
is a solution of the SDE
$$
X_y=2\,\int_0^y\sqrt{X_s}\,d\beta_s+2(\,y\wedge 1\,),
$$
in other words,
\begin{enumerate}
\item $\{L^{1-y}_{H_1(B)}(B):\ 0\leq y\leq 1\}$ is a 2-dimensional
squared Bessel process starting from 0,
\item   $\{L^{1-y}_{H_1(B)}(B):\ y\geq 1\}$ is a 0-dimensional
squared Bessel process with the starting value $L^{0}_{H_1}(B)$
obtained from (i).
\end{enumerate}
}
\end{theorem}

\begin{theorem}
\label{r-k2}
{\sl Assume that $B^{(\mu)}_0=0$ and let $L^y_t(B^{(\mu)})$ be the
local time of $B^{(\mu)}$ at level $y$ (with
respect to the Lebesgue measure) up to time $t.$ Define 
the total local time of $B^{(\mu)}$ at level $y$
via
$$
L^y_\infty(B^{(\mu)}):=\lim_{t\to\infty}L^y_t(B^{(\mu)}).
$$
Then
\begin{equation}
\label{rk1}
\{L^{-y}_\infty(B^{(\mu)}):\ y\geq 0 \}\
{\mathop=^{\rm{(d)}}}\
\{Z^{(0,2\mu)}_y:\ y\geq 0 \},
\end{equation}
and
    \begin{equation}
    \label{rk2}
    \{L^y_\infty(B^{(\mu)}):\ y\geq 0 \}\
    {\mathop=^{\rm{(d)}}}\
    \{Z^{(2,2\mu)}_y:\ y\geq 0\},
    \end{equation}
    where $Z^{(\delta,2\mu)}, \delta=0,2,$ are solutions of the SDE
$$
    dX_t= 2\sqrt{X_t}\,dB_t+(\delta -2\mu X_t)\,dt,
$$
    respectively, with the initial value $X_0$
    exponentially distributed with parameter $\mu.$
    In fact, the identities (\ref{rk1}) and (\ref{rk2}) hold jointly,
    with $Z^{(0,2\mu)}_0=Z^{(2,2\mu)}_0$ but
otherwise $Z^{(0,2\mu)}$ and $Z^{(2,2\mu)}$ are independent.
}
\end{theorem}

Our final Ray--Knight theorem is for Bessel processes.
The first part is
formulated only for 3-dimensional Bessel
process, and in the second part we take the dimension parameter $\delta>2.$
Let $L^y_\infty(R^{(\delta)})$ denote the total local time at $y$ of
the Bessel process $R^{(\delta)}$ (taken with respect to the Lebesgue measure).
\begin{theorem}
\label{r-k3}
{\bf a:}\quad{\sl  Assume that $R^{(3)}$ is started at 0. Then
$$
\{L^y_\infty(R^{(3)}):\ y\geq 0 \}\
{\mathop=^{\rm{(d)}}}\
\{Z^{(2)}_y:\ y\geq 0 \},
$$
where $Z^{(2)}$ denotes the squared Bessel process of dimension 2,
started from 0, i.e., $Z^{(2)}$ satisfies the SDE
$$
dX_y= 2\sqrt{X_y}\,dB_y+2\,dy.
$$
\hfill\break\hfill
{\bf b:}\quad Assume that $\delta>2$ and $R^{(\delta)}_0=0.$ Then
$$
\{L^y_{H_1}(R^{(\delta)}):\ 0\leq y\leq 1 \}\
{\mathop=^{\rm{(d)}}}\
\Big\{\frac{1}{(\delta-2)y^{\delta-3}}\ {\widehat Z^{(2)}_{y^{\delta-2}}}
:\ 0\leq y\leq 1 \Big\}
$$
where $\widehat Z^{(2)}$ denotes the 2-dimensional
squared Bessel bridge (from 0 to
0 and of length 1).
}
\end{theorem}

\bibliographystyle{plain}
\bibliography{yor1}
\end{document}